\documentclass[11pt]{article}

\usepackage{amssymb,graphics,amsmath,amsthm,amsopn,amstext,amsfonts,color,fullpage,fancybox,bm}
\usepackage{subfig,placeins} 
\usepackage{graphicx,color,float,epstopdf}
\newcommand{\gap}{\vspace{0.1in}}
\newcommand{\epc}{\hspace{1pc}}
\newcommand{\thalf}{{\textstyle{\frac{1}{2}}}}
\newcommand{\wh}{\widehat}
\newcommand{\wt}{\widetilde}
\newcommand{\onebld}{\mathbf{1}}

\newcommand{\E}{{\rm I\!E}}

\newcommand{\IP}{{\rm I\!P}}
\newcommand{\Id}{{\rm I\!I}}

\newtheorem{theorem}{Theorem}
\newtheorem{proposition}[theorem]{Proposition}
\newtheorem{lemma}[theorem]{Lemma}
\newtheorem{corollary}[theorem]{Corollary}

\title{On the Pervasiveness of Difference-Convexity in \\
Optimization and Statistics\footnote{
All authors are affiliated with the Daniel J.\ Epstein Department of Industrial and Systems
Engineering, University of Southern California, Los Angeles, California 90089-0193, U.S.A.
{\tt Emails: nouiehed@usc.edu; jongship@usc.edu; razaviya@usc.edu}}}

\date{Original: February 2017; 
Last revision: April 2018}

\author{Maher Nouiehed \and Jong-Shi Pang\footnote{This work of this author was based on
research partially supported by the U.S.\ National Science Foundation grants CMMI 1538605 and IIS-1632971.}
\and Meisam Razaviyayn}

\begin{document}

\maketitle

\begin{abstract}
\noindent With the increasing interest in applying the methodology of difference-of-convex (dc)
optimization to diverse
problems in engineering and statistics, this paper establishes the dc property of many
functions in various areas of applications not previously known to be of this class.  Motivated
by a quadratic programming based
recourse function in two-stage stochastic programming, we show that the (optimal) value function
of a copositive (thus not necessarily convex) quadratic program is dc on the domain of finiteness
of the program when the matrix in the objective function's quadratic term and the constraint matrix
are fixed.  The proof of this result is based on a dc decomposition of a piecewise LC$^1$ function
(i.e., functions with Lipschitz gradients).  Armed with these new results and known properties of
dc functions existed in the literature, we show that many composite statistical functions in risk analysis,
including the value-at-risk (VaR), conditional value-at-risk (CVaR), Optimized Certainty Equivalent
(OCE), and the expectation-based, VaR-based, and CVaR-based
random deviation functionals are all dc.  Adding the known class of dc surrogate sparsity functions that are employed
as approximations of the $\ell_0$ function in statistical learning, our work significantly expands the classes
of dc functions and positions them for fruitful applications.
\end{abstract}

\section{Introduction}


Long before their entry into the field of optimization in the early 1980's
\cite{HThoai99,Hiriart-Urruty85,PhamLeThi97,Tuy87}, difference-of-convex (dc) functions have been studied
extensively in the mathematics literature; the 1959 paper \cite{Hartman59}
cited a 1950 paper \cite{Alexandroff50} where dc functions were considered.  The paper \cite{Hartman59}
contains a wealth of fundamental results on dc functions that lay the foundation for this class of non-convex
functions.  While focused on the more general class of ``delta-convex functions'' in abstract spaces,
the thesis \cite{VeselyZajicek89} contains the very important mixing property of dc functions that in today's
language is directly relevant to the dc property of {\sl piecewise functions}.  The most recent
paper \cite{BBorwein11} adds to this
literature of the mathematics of dc functions with a summary of many existing properties of dc functions.
As noted in the last paper, the mapping that Nash employed to show
the existence of a mixed equilibrium strategy in his celebrated 1951 paper \cite{Nash51} turns out to be defined by
dc functions.  This provides another evidence of the relevance of dc functions more than half a century ago.
In the optimization literature, applications of the dc methodology to nonconvex optimization problems are well
documented in the survey papers \cite{LeThiPham05,LeThiPham14} and in scattered papers by the pair of authors
of the latter papers and their collaborators; adding to these surveys, the paper \cite{JaraMoroniWaechterPang16}
discusses an application of dc programming to
the class of linear programs with complementarity constraints;
the most recent paper \cite{LeThiPham16} documents many contemporary applied problems in diverse engineering
and other disciplines.

\gap

Our own interest in dc functions stemmed from the optimization of some physical layer problems
in signal processing and communication \cite{AhnPangXin17,AlvaradoScutariPang14,PangRazaviyaynAlvarado16}.
Most worthy of note in these references are the following.  In \cite{PangRazaviyaynAlvarado16}, a novel class of
dc functions was identified and an iterative algorithm was described to compute a directional stationary point
of a convex constrained dc program; extensions of the algorithm to dc constraints were also presented.
In \cite{AhnPangXin17},
a unified dc representation was given for a host of surrogate sparsity functions that were employed as
approximations of the $\ell_0$ function in statistical learning; such a representation further confirms the fundamental
importance of dc functions in the latter subject that is central to today's field of big-data science
and engineering.  The paper \cite{PangTao17} investigates decomposition methods solving a class of multi-block
optimization problems with coupled constraints and partial dc-structure.

\gap

The present paper was initially motivated by
the desire to understand the dc property of composite risk/deviation functions arising from financial engineering
\cite{RUryasev13,RUZabarankin06a,RUZabarankin06b}; see
Section~\ref{sec:composite risk fncs} for a formal definition of these functions.  Roughly speaking,
these are functions defined as the compositions of some well-known statistical
quantities, such as variance, standard deviation, and quantiles, with some random functionals such as the uncertain
return of an investment portfolio \cite{RUryasev00,RUryasev02,SarykalinUryasev08} or the second-stage recourse function
\cite{PangSenShanbhag16} in stochastic programming \cite{BirgeLouveaux97,ShapiroDentchenvaRuszczynski09}.
Initially, we were intrigued by the question of whether the value-at-risk (VaR) functional of a random portfolio return
was a dc function of the asset holdings.  It turns out that this question already has an affirmative answer given
in the paper \cite{Wozabal12} where a formula
linking the VaR and CVaR is obtained.  In our work, we provide an alternative expression connecting these two risk quantities that
is based on linear programming duality.
We also extend our formula to the more general context of the optimized certainty equivalent introduced by Ben Tal and
Teboulle \cite{BenTalTeboulle86,BenTalTeboulle87} that predates the work of Rockafellar and Uryasev \cite{RUryasev00,RUryasev02}.
Our next
investigation pertains to a composite risk function involving a recourse function defined by a random quadratic program
parameterized by the first-stage decision variable.  Our analysis pertains to the general problem where the latter quadratic
recourse program is nonconvex.  A second contribution of our work is a detailed proof showing that the (optimal)
value function of a quadratic program (QP) is dc on the domain of finiteness
of the program when the matrix in the objective function's quadratic term and the constraint matrix
are fixed, under the assumption that the former matrix is copositive on the recession cone of the constraint region.  Such
a copositive assumption is essential because without it the quadratic program is unbounded below on any non-empty feasible set.
In turn, the proof of the said dc property of the QP value function is based on an explicit dc decomposition of a piecewise
function with Lipschitz gradients that is new by itself.  In addition to these specialized results pertaining to statistical
optimization and two-stage
stochastic programming, we obtain a few general dc results that supplement various known facts in the literature; e.g.,
the decompositions in part (c) of Lemma~\ref{lm:general dc properties} and Propositions~\ref{pr:composite incr cvx}
and \ref{pr:composite log}.

\subsection{Significance of the dc property}

Besides the mathematical interest, the dc property of a function can be used profitably for the design of convex program based
optimization algorithms.  Indeed, the backbone of the classical dc algorithm \cite{LeThiPham05,LeThiPham14,PhamLeThi97} and its recent
enhancements \cite{PangRazaviyaynAlvarado16,PangTao17} is a given dc decomposition of the functions involved in the
optimization problem.  Such a decomposition provides a convenient convex majorization of the functions that can be used as
their surrogates to be optimized \cite{Razaviyayn14,RazaviyaynHongLuo13}.  In a nutshell, the benefit of a dc-based
iterative algorithm is that it provides a descent algorithm without either a line search or a trust region step; as a result,
parallel and/or distributed implementations of such an algorithm can be easily designed without centralized
coordination \cite{RazaviyaynHongLuoPang16,RazaviyaynHongLuoPang14,ScutariFacchineiPalomarPangSong14,ScutariAlvaradoPang14}
when the given problem has certain partitioned structure.
More interestingly, a certain class of dc programs is the only class of nonconvex, nondifferentiable programs
for which a directional derivative based stationary point can be computed \cite{PangRazaviyaynAlvarado16};
such a stationarity concept is the sharpest one among all ``first-order'' stationarity concepts.
In particular, directional stationarity is significantly sharper than the convex-analysis based
concept of a critical point of a dc program \cite{LeThiPham05,LeThiPham14,PhamLeThi97}, which is a relaxation of a Clarke
stationary point \cite{Clarke83} under the dc property.   Here
sharpness refers to the property that under various first-order stationarity definitions, the
corresponding sets of stationary points contain the set of directional stationary points.
The papers \cite{AhnPangXin17,ChangHongPang17} examine how directional stationarity is instrumental in the characterization of local
minimizers of nonconvex and nondifferentiable optimization problems under piecewise linearity and/or second-order conditions.

\gap

The dc algorithm and its variants have been employed in many applied contexts; see the recent survey \cite{LeThiPham16}. For many
existing applications, the resulting convex programs can be solved very easily; see e.g.\ \cite{GotohTakedaKono17,LeThiPhamVo15} in the
area of sparse optimization that has attracted much interest in recent years.  In general, deciding whether a given nonconvex
function is dc is not necessarily an easy task.  In situations
where a function can be shown to be dc,
a dc decomposition offers the first step to design an efficient algorithm by readily providing a convex majorant of
the function to be minimized and enabling the investigation and computation of sharp stationary solutions.
A novel case in point is the family of nonconvex, nondifferentiable composite programs arising from piecewise affine statistical regression
and multi-layer neural networks with piecewise affine activation functions.  For these applied problems, the convex
subprograms are not straightforward to be solved; nevertheless, they are amenable to efficient solution by
a semi-smooth Newton method \cite[Section~7.5]{FacchineiPang03}.  Details of this dc approach to the numerical
solution of these advanced statistical learning problems can be found in the most recent reference \cite{CuiPangSen18}.

\gap

In summary, a dc decomposition offers both a computational venue and a theoretical framework for the understanding and numerical solution
of nonconvex nondifferentiable optimization problems.
In this paper we collect in one place many nonconvex functions that are not previously known to be of the dc kind and establish their dc property.
These function arise from 3 different areas: (a) composite risk functionals, (b) statistical estimation and learning,
and (c) quadratic recourse in stochastic programming.  Further study of how the dc property of these functions can benefit algorithmic design for solving the
optimization problems involving these functions is beyond the scope of this paper.  As word of caution: the obtained dc decompositions for the functions
studied in this paper may not be the most conducive for numerical use; nevertheless, they offer a formal demonstration of the dc property of the
functions that are not previously known to have this property.

\subsection{Organization of the paper}

In the next section, we introduce several classes of composite risk and statistical functions to be studied subsequently.
These include the renowned Value-at-Risk (VaR) and Conditional Value-at-Risk (CVaR) \cite{RUryasev00,RUryasev02}, their extensions
to a utility based Optimized Certainty Equivalent (OCE) \cite{BenTalTeboulle07,BenTalTeboulle87,BenTalTeboulle86},
and maximum likelihood functions
(Subsection~\ref{subsec:one-parameter exponential}) derived from
one-parameter exponential densities composite with a statistical estimation model. 
In Section~\ref{sec:dc random functionals}, we give proofs of the difference-convexity of the composite risk and
statistical functions, providing in particular an alternative expression connecting the
VaR with the CVaR in Subsection~\ref{subsec:CVaR and VaR} that is extended to the OCE in Subsection~\ref{subsec:OCE}.
Subsection~\ref{subsec:statistical functions} deals with the composite statistical functions such as the variance,
standard deviations and certain composite density functions.
Section~\ref{sec:quadratic value fnc} shows that the optimal
objective value of a copositive quadratic program for a fixed constraint matrix is a dc function of the linear term of the objective
and the right-hand side of the constraints; such a value function arises as the recourse of a two-stage stochastic program with
the copositivity assumption generalizing the positive semi-definiteness of the quadratic form and making the objective function nonconvex
in general.  Specialized
to a linear program, this dc property is a new result in the vast literature of linear programming theory and provides a first
step in expanding the two-stage stochastic programming domain beyond the much focused paradigm of linear recourse.
Lastly in Section~\ref{sec:folded concave} we give a necessary and sufficient condition for a univariate folded concave function to be dc.
Such a function provides a unification of all the surrogate sparsity functions that have a fundamental role to play in nonconvex sparsity
optimization.

\section{Composite Risk and Statistical Functions} \label{sec:composite risk fncs}

%

For a given scalar $\lambda > 0$ and decision variable $x \in \mathbb{R}^n$, consider the function
\[
{\cal R}_{\lambda}(x) \, \triangleq \, \E {\cal Z}(x;\wt{\omega}) +
\lambda \, {\cal D}\left[ {\cal Z}(x;\wt{\omega}) \right]
\]
where $\E$ is the expectation operator with respect to the random variable $\wt{\omega}$ defined on the
probability space $( \Omega,{\cal F}, \IP )$, with $\Omega$ being the sample space,
${\cal F}$ being the $\sigma$-algebra generated by subsets of $\Omega$, and $\IP$ being a probability measure
defined on ${\cal F}$; $\lambda > 0$ is a given parameter that
balances the expectation (for risk neutrality) and the deviation measure ${\cal D}$ (representing risk aversion);
for a random variable ${\cal Z}$, ${\cal D}( {\cal Z} )$ is a expectation-based, CVaR-based, or VaR-based deviation measure:

\gap



\begin{tabular}{lll}
& \hspace{0.1in} \textbf{Expectation Based} & \hspace{0.1in} \textbf{(C)VaR Based} \\ [5pt]
$\bullet $ Variance:  & \hspace{0.1in}
$\sigma^2({\cal Z}) \triangleq \E\left[ \, {\cal Z} - \E{\cal Z} \, \right]^2$  
& \hspace{0.1in}
$\E\left[ \, {\cal Z} - \mbox{(C)VaR}({\cal Z}) \, \right]^2$  \\ [5pt]
$\bullet $ Standard Deviation: & \hspace{0.1in}
$\sigma({\cal Z}) \triangleq \sqrt{\sigma^2({\cal Z})}$ 
& \hspace{0.1in}
$\sqrt{\E\left[ \, {\cal Z} - \mbox{(C)VaR}({\cal Z}) \, \right]^2}$  \\ [5pt]
$\bullet $ Absolute Semi-Deviation (ASD):  & \hspace{0.1in}
$\mbox{ASD}({\cal Z}) \triangleq \E \left[ {\cal Z} - \E{\cal Z} \right]_+$ 
& \hspace{0.1in}
$\E \left[ {\cal Z} - \mbox{(C)VaR}({\cal Z}) \right]_+$  \\ [5pt]
$\bullet $ Absolute Deviation (AD): & \hspace{0.1in}
$\mbox{AD}({\cal Z}) \triangleq \E \left| \, {\cal Z} - \E{\cal Z} \, \right|$ 
& \hspace{0.1in}
$\E \left| \, {\cal Z} - \mbox{(C)VaR}({\cal Z}) \, \right|$
\end{tabular}














\gap

where $[ \bullet ]_+ \triangleq \max( 0, \bullet )$, and the expressions of $\mbox{CVaR}_{\alpha}(\bullet )$
and $\mbox{VaR}_{\alpha}(\bullet )$ are as follows:
\[ \begin{array}{lll}
\mbox{CVaR}_{\alpha}({\cal Z}) & = & \displaystyle{
\operatornamewithlimits{\mbox{minimum}}_{t \in \mathbb{R}}
} \, \left[ \, t + \displaystyle{
\frac{1}{1 - \alpha}
} \, \E\left[ \, {\cal Z} - t \, \right]_+ \, \right] \\ [0.2in]
\mbox{VaR}_{\alpha}({\cal Z}) & = & \mbox{minimum}\left\{ \, t^{\, \prime} \, \mid \, t^{\, \prime} \, \in \displaystyle{
\operatornamewithlimits{\mbox{argmin}}_{t \in \mathbb{R}}
} \, \left[ \, t + \displaystyle{
\frac{1}{1 - \alpha}
} \, \E\left[ \, {\cal Z} - t \, \right]_+ \, \right] \, \right\}.
\end{array} \]
The case where ${\cal D}({\cal Z})$ is itself $\mbox{CVaR}_{\alpha}({\cal Z})$
or $\mbox{VaR}_{\alpha}({\cal Z})$ is also covered by our analysis.
Unlike the equality $\mbox{AD}({\cal Z}) = 2\mbox{ASD}({\cal Z})$,
it is in general not true that
$\E \left| {\cal Z} - \mbox{CVaR}_{\alpha}({\cal Z}) \right| = 2 \E \left[ {\cal Z} - \mbox{CVaR}_{\alpha}({\cal Z}) \right]_+$.
For a recent survey on the connection between
risk functions and deviation measures and their applications in risk management and statistical estimation, see \cite{RUryasev13}.

\gap

Admittedly, the (C)Var-based deviation
measure is a non-traditional quantity that is not commonly employed in risk analysis.  Nevertheless, as a quantile of a random
variable that includes for instance the median as a special case, we feel that it is important to consider the deviation from
such a fundamental statistical
quantity and understand its generalized convexity properties (if applicable) when its optimization is called for.
The resulting deviations
are the analogs of the classical variance family of deviations that are based on the mean of the random variable.

\gap

Extending the (C)VaR, the Optimized Certainty Equivalent (OCE) of a random variable is defined
by a proper, concave,
non-decreasing, lower semi-continuous utility function
$u : \mathbb{R} \to [ \, -\infty, \, \infty \, )$ with a non-empty effective domain
$\mbox{dom}(u) \triangleq \left\{ \, t \in \mathbb{R} \mid u(t) > -\infty \, \right\}$
such that $u(0) = 0$ and $1 \in \partial u(0)$, where $\partial u$ denotes the subdifferential map of $u$.  Thus in particular,
\[ \begin{array}{ll}
u(t) \, \geq \, 0, & \forall \, t \, \geq \, 0, \epc \mbox{and} \\ [5pt]
u(t) \, \leq \, t, & \forall \, t \, \in \, \mathbb{R}.
\end{array}
\]
Let ${\cal U}$ be the family of these univariate utility functions.
For an essentially bounded random variable ${\cal Z}$ satisfying $\displaystyle{
\sup_{\omega \in \Omega}
} \, | \, {\cal Z}(\omega) \, | \, < \, \infty$, the {\sl optimized certainty equivalent} (OCE) of ${\cal Z}$
with respect to a utility function $u \in {\cal U}$ is defined as
\[
{\cal O}_u({\cal Z}) \, \triangleq \, \displaystyle{
\sup_{\eta \in \mathbb{R}}
} \, \left[ \, \eta + \E u( {\cal Z} - \eta ) \, \right].
\]
The choice of $u(t) \triangleq \displaystyle{
\frac{1}{1 - \alpha}
} \, \min(0,t)$ yields ${\cal O}_u({\cal Z}) = -\mbox{CVaR}_{\alpha}(-{\cal Z})$.  Proposition~2.1 in \cite{BenTalTeboulle07}
shows that for a random variable ${\cal Z}$ whose support is a compact interval, then the supremum in ${\cal O}_u({\cal Z})$ is attained;
in this case, we may consider the largest such maximizer,
\[
m_u({\cal Z}) \, \triangleq \, \max\left\{ \, \eta^{\, \prime} \, \mid \, \eta^{\, \prime} \, \in \displaystyle{
\operatornamewithlimits{\mbox{argmax}}_{\eta \in \mathbb{R}}
} \, \left[ \, \eta + \E u( {\cal Z} - \eta ) \, \right] \, \right\},
\]
which is the utility-based extension of VaR that pertains to a given choice of $u$.  As explained in \cite{BenTalTeboulle07},
the deterministic quantity $m_u({\cal Z})$ can be interpreted as the largest optimal allocation between present and future consumption
if ${\cal Z}$ represents an uncertain income of ${\cal Z}$ dollars, largest as a way to break ties among multiple optimal allocations
if such allocation is not unique.

\subsection{Quadratic recourse function} \label{subsec:quadratic recourse}

Besides the bilinear ${\cal Z}(x;\omega) = x^T\omega$ that is quite common in portfolio management with $\omega$ representing
the uncertain asset returns and $x$ the holdings of the assets, we shall treat carefully a quadratic recourse function given by
\begin{equation} \label{eq:QP recourse-definition}
\begin{array}{llll}
\psi(x;\omega) & \triangleq & \displaystyle{
\operatornamewithlimits{\mbox{minimum}}_{z}
} & \left[ \, f(\omega) + G(\omega)x \, \right]^Tz + \thalf \, z^TQz\\ [7pt]
& & \mbox{subject to} & C(\omega)x + Dz \, \geq \, \xi(\omega),
\end{array} \end{equation}
where $Q \in \mathbb{R}^{m \times m}$ is a symmetric, albeit not necessarily positive or negative semi-definite, matrix;
$D$ is a $k \times m$ matrix, $f : \Omega \to \mathbb{R}^m$ and $\xi : \Omega \to \mathbb{R}^k$ are vector-valued random functions,
and $G : \Omega \to \mathbb{R}^{m \times n}$ and $C : \Omega \to \mathbb{R}^{k \times n}$ are matrix-valued random functions.
Besides the main dc result of the composite function ${\cal R}_{\lambda}(x)$,
the proof that the value function $\psi(\bullet;\omega)$ is dc on its domain of finiteness is a major contribution
of this work that is of independent interest.  This result will be discussed in detail in Section~\ref{sec:quadratic value fnc}.
There are several noteworthy points of our analysis: (a) the matrix $Q$ is not required to be
positive semi-definite; thus we allow our recourse function to be derived from an indefinite quadratic program; (b) the
first-stage variable $x$ appears in both the objective function and the constraint; this is distinguished from much of the
stochastic (linear) programming literature where $x$ appears only in the constraint; and (c) it follows from our result that
the value function of a linear program:
\[ \begin{array}{cccl}
\varphi(b,c) & \triangleq & \displaystyle{
\operatornamewithlimits{\mbox{minimum}}_{z}
} & c^Tz \\ [5pt]
& & \mbox{subject to} & Az \, = \, b \epc \mbox{and} \epc z \geq 0,
\end{array} \]
is a dc function on its domain of finiteness.  This is a new result by itself because existing results in parametric linear programming
deal only with the concavity/convexity of $\varphi(b,c)$ when $b$ ($c$ respectively) is fixed; there does not exist a
(non-)convexity analysis of the optimal objective value as a function jointly of $b$ and $c$.

\gap

Like the deviations from the (C)VaR, recourse-function based deviations are not common in the stochastic programming literature.
Part of the reason for this lack of attention might be due to the computational challenge of dealing with the recourse function itself,
which is further complicated when coupled with the statistical functions, such as variances, standard deviations, or semi-deviations.
Hopefully, understanding the structural properties of the composite recourse-based deviations could open a path for solving the advanced
stochastic programs for risk-averse players who might be interested in reducing their risk exposure of deviation from the second-stage decisions.
In such a situation, the composite deviations would provided a reasonable measure of the risk to be reduced.

\gap

In the analysis of the value function $\psi(\bullet;\omega)$,
we are led to a detailed study of the dc property of piecewise functions.  Specifically, a continuous
function $\theta$ defined on an open set ${\cal O} \subseteq \mathbb{R}^n$ is piecewise C$^k$ \cite[Definition~4.5.1]{FacchineiPang03}
for an integer $k \geq 0$ if there exist finitely many C$^k$ functions $\{ \theta_i \}_{i=1}^I$ for some integer $I > 0$,
all defined on ${\cal O}$, such that $\theta(x) \in \{ \theta_i(x) \}_{i=1}^I$ for all $x \in {\cal O}$.  A major focus of our
work is the case of a piecewise quadratic (PQ) $\theta$, which has each $\theta_i$ being a (possibly nonconvex) quadratic function.
It is a well-known fact that a general quadratic function must be dc; subsequently, we extend this fact to a piecewise quadratic function
with an explicit dc representation in terms of the pieces.  While our analysis relies on several
basic results of dc functions that can be found in \cite{BBorwein11,Hartman59} (see also \cite[Chapter~4]{Tuy16}),
we also discover a number of new results
concerning PQ functions that are of independent interest.  Summarized below,
these results are the pre-requisites to establish the dc property of the composite-risk functions with quadratic recourse.


\gap

$\bullet $ The optimal objective value of the quadratic program (QP):
\begin{equation} \label{eq:copositive qp}
\begin{array}{lll}	
\mbox{qp}_{\rm opt}(q,b) \, \triangleq & \displaystyle{
\operatornamewithlimits{\mbox{minimum}}_{z}
} & \zeta(z) \, \triangleq \, q^Tz + \thalf \, z^TQz\\ [7pt]
& \mbox{subject to} & z \, \in \, {\cal P}_D(b) \, \triangleq \, \left\{ \, z \, \in \, \mathbb{R}^m \, \mid \, Dz \, \geq \, b \, \right\},
\end{array} \end{equation}
is a dc function of $(q,b)$ on the domain
$\mbox{dom}(Q,D) \triangleq \left\{ ( q,b ) \in \mathbb{R}^{m+k} \mid -\infty < \mbox{qp}_{\rm opt}(q,b) < \infty \right\}$
of finiteness of the problem, for a fixed pair $(Q,D)$ with $Q$ being a symmetric matrix that is {\sl copositive} on
the recession cone of the feasible region ${\cal P}_D(b)$; i.e., provided that $v^TQv \geq 0$ for all
$v \in D_{\infty} \triangleq \left\{ v \in \mathbb{R}^m \mid Dv \geq 0 \right\}$.  We let $\mbox{qp}_{\rm sol}(q,b)$ denote
the optimal solution set of (\ref{eq:copositive qp}), which is empty for $(q,b) \not\in \mbox{dom}(Q,D)$.  It turns out
that the analysis of (\ref{eq:copositive qp}) in such a copositive case is not straightforward and uses significant background
about the problem and polyhedral theory.

\gap

$\bullet $ Motivated by the value function $\mbox{qp}_{\rm opt}(q,b)$,
we obtain an explicit min-max (dc) representation of a general (not necessarily convex)
piecewise quadratic function with given pieces, extending the work
\cite{Sun92} that studies the special case when such a function is convex and also the
max-min representation of a piecewise linear function \cite{Ovchinnikov02,Scholtes12}, as well as
the so-called {\sl mixing property} of dc functions \cite[Lemma~4.8]{VeselyZajicek89}
on open convex sets extended to the family of PQ functions whose domains are closed sets.

\subsection{One-parameter exponential densities} \label{subsec:one-parameter exponential}

The discussion of the following statistical modeling is drawn from \cite{Sen17}.  A random variable $Y$ belongs
to a {\sl one-parameter exponential family} if its
density (or mass) function can be written in the form
\[
g(y;\theta) \, = \, a(y) \exp\left\{ \, y \theta - b(\theta) \, \right\},
\]
where $\theta$ is the canonical parameter, and $a$ and $b$ are given functions with
$b$ being convex and increasing.  In the presence of covariates $X$, we model the conditional
distribution of $Y|X = x$ as a one-parameter exponential family member with parameter $\theta = m(x;\Theta)$
that depends on the realization $x$ of $Y|X$ and where $m(\bullet;\Theta)$ is a parametric statistical (e.g.\ linear)
model with the parameter  $\Theta \in \mathbb{R}^n$ in
the latter model being computed by maximzing the expected log-likelihood function
$\E_{Y|X}\left[ \, \log g(y;m(x;\Theta)) \, \right]$.
When discretized with respect
to the given data $\{ (y_s;x^s) \}_{s=1}^N$, the latter optimization problem is equivalent to
\[
\displaystyle{
\operatornamewithlimits{\mbox{maximize}}_{\Theta}
} \, \displaystyle{
\frac{1}{N}
} \displaystyle{
\sum_{s=1}^N
} \, \left[ \ y_s \, m(x^s;\Theta) - b \circ m(x^s;\Theta) \, \right].
\]
In the recent paper \cite{HahnBanergjeeSen16}, a piecewise affine statistical estimation model was proposed with
$m(\bullet;\Theta)$ being a piecewise affine function.  Per the representation results of \cite{Scholtes12,Sun92},
every such function can be expressed as the difference of two convex piecewise affine functions, and is thus dc.
This motivates us to ask the question of whether the composite function $b \circ m(\bullet;\Theta)$ is dc, and
more generally, if
$\E\log f(\Theta;\wt{\omega})$ is dc if $f(\bullet;\omega)$ is dc for each $\omega \in \Omega$.

\section{Proof of Difference-Convexity}  \label{sec:dc random functionals}

We show that if the random function ${\cal Z}(x;\omega) = p(x,\omega) - q(x,\omega)$ where
$p(\bullet,\omega)$ and $q(\bullet,\omega)$ are both convex functions on a domain ${\cal D}$ for every fixed $\omega \in \Omega$,
then all the expectation based and (C)VaR based risk measures ${\cal R}_{\lambda}(x)$
are dc on ${\cal D}$, so are the OCE extensions ${\cal O}_u({\cal Z}(x;\wt{\omega}))$
and $m_u({\cal Z}(x;\wt{\omega}))$ with a piecewise linear utility function $u$.   In turn,
it suffices to show that the following functions are dc with ${\cal Z} = {\cal Z}(x;\wt{\omega})$:

\gap

$\bullet $ $\mbox{CVaR}_{\alpha}({\cal Z})$ and $\mbox{VaR}_{\alpha}({\cal Z})$;

\gap

$\bullet $ the OCE extensions ${\cal O}_u({\cal Z})$
and $m_u({\cal Z})$ with a piecewise linear utility function $u$;

\gap

$\bullet $ $\sigma^2({\cal Z})$ and $\sigma({\cal Z})$.

\gap

Once these are shown, using known properties of
dc functions (such as the nonnegative part of a dc function is dc), we can readily establish the dc-property of the
following functions and also yield their dc decompositions:

\gap

$\bullet $ $\mbox{ASD}({\cal Z}) = \thalf \mbox{AD}({\cal Z})$;

\gap

$\bullet $ $\E\left[ {\cal Z} - \mbox{(C)VaR}_{\alpha}({\cal Z}) \right]^2$
and $\sqrt{\E\left[ {\cal Z} - \mbox{(C)VaR}_{\alpha}({\cal Z}) \right]^2}$; and

\gap

$\bullet $
$\E\left[ {\cal Z} - \mbox{(C)VaR}_{\alpha}({\cal Z}) \right]_+$ and
$\E \left| {\cal Z} - \mbox{(C)VaR}_{\alpha}({\cal Z}) \right|$.

\gap

Among the former three families of composite functions, the proof of
$\mbox{VaR}_{\alpha}({\cal Z}(x;\wt{\omega}))$, $m_u({\cal Z}(x;\wt{\omega}))$, $\sigma({\cal Z}(x;\wt{\omega}))$
requires the random variable $\wt{\omega}$  to be discretely distributed.

\subsection{CVaR and VaR} \label{subsec:CVaR and VaR}

The following result shows that $\mbox{CVaR}_{\alpha}(f(x,\wt{\omega}))$ is a dc function if $f(\cdot,\omega)$ is dc
for fixed realization $\omega$.

\begin{proposition} \label{pr:CVaR} \rm
For every $\omega \in \Omega$, let $f(x,\wt{\omega}) = p(x,\wt{\omega}) - q(x,\wt{\omega})$ be the dc decomposition of $f(\bullet,\omega)$
on a convex set ${\cal D} \subseteq \mathbb{R}^n$.  Then for every $\alpha \in ( 0,1 )$,
\[
\mbox{CVaR}_{\alpha}(f(x,\wt{\omega})) \, = \, \underbrace{\displaystyle{
\min_{t \in \mathbb{R}}
} \, \left\{ \, t + \displaystyle{
\frac{1}{1 - \alpha}
} \, \E\max\left( \, p(x,\wt{\omega}) - t, q(x,\wt{\omega}) \, \right) \, \right\}}_{\mbox{cvx in $x$}} -
\underbrace{\displaystyle{
\frac{1}{1 - \alpha}
} \, \E q(x,\wt{\omega})}_{\mbox{cvx in $x$}}.
\]
\end{proposition}

\noindent {\bf Proof.} The equality is fairly straightforward.  The convexity of the minimum is due to the joint convexity
of the function $(x,t) \mapsto t + \displaystyle{
\frac{1}{1 - \alpha}
} \, \E \max\left( \, p(x,\wt{\omega}) - t, q(x,\wt{\omega}) \, \right)$.  \hfill $\Box$

\gap

A formula that connects $\mbox{VaR}_{\alpha}({\cal Z})$ to $\mbox{CVaR}_{\alpha}({\cal Z})$ for a discretely
distributed random variable ${\cal Z}$ was obtained in \cite{Wozabal12}; this formula can be used to establish
the dc property of $\mbox{VaR}_{\alpha}(f(x,\wt{\omega}))$ when $f(\bullet,\wt{\omega})$ is dc..
In what follows, we derive an alternative expression connecting $\mbox{VaR}({\cal Z})$ and $\mbox{CVaR}({\cal Z})$ using simple
linear programming duality when the sample space
$\Omega = \left\{ \omega^1, \cdots, \omega^S \right\}$ for some integer $S > 0$.
Let $\left\{ p_1, \cdots, p_S \right\}$ be the associated family of probabilities of the discrete realizations
of the random variable $\wt{\omega}$.  Writing $z_s \triangleq f(x,\omega^s)$, we then have
\[ \begin{array}{l}
\mbox{VaR}_{\alpha}( f(x,\wt{\omega}) ) \, = \, \min\left\{ \, t^{\, \prime} \, \mid \, t^{\, \prime} \, \in \displaystyle{
\operatornamewithlimits{\mbox{argmin}}_{t \in \mathbb{R}}
} \, \left( \, t + \displaystyle{
\frac{1}{1 - \alpha}
} \, \displaystyle{
\sum_{s=1}^S
} \ p_s \, \left[ \, z_s - t \, \right]_+ \, \right) \, \right\} \\ [0.2in]
= \, \displaystyle{
\operatornamewithlimits{\mbox{minimum}}_{t,u}
} \hspace{0.3in} t \hspace{2in} \mbox{(by representation of the argmin)} \\ [3pt]
\epc \mbox{subject to} \epc t + \displaystyle{
\frac{1}{1 - \alpha}
} \, \displaystyle{
\sum_{s=1}^S
} \, p_s \, u_s \, \leq \, \mbox{CVaR}_{\alpha}( f(x,\wt{\omega}) ) \\ [0.2in]
\hspace{0.85in} \left. \begin{array}{rcl}
t + u_s & \geq & z_s \\ [3pt]
u_s & \geq & 0
\end{array} \right\} \epc s \, = \, 1, \cdots, S \\ [0.2in]
= \, \displaystyle{
\operatornamewithlimits{\mbox{maximum}}_v
} \hspace{0.15in} -v_0 \, \mbox{CVaR}_{\alpha}( f(x,\wt{\omega}) ) + \displaystyle{
\sum_{s=1}^S
} \, v_s \, z_s \epc \mbox{(by linear programming duality)} \\ [0.2in]
\epc \mbox{subject to} \epc -v_0 + \displaystyle{
\sum_{s=1}^S
} \, v_s \, = \, 1 \\ [0.2in]
\hspace{1in} \displaystyle{
\frac{-p_s}{1 - \alpha}
} \, v_0 + v_s \, \leq \, 0, \epc s \, = \, 1, \cdots, S \\ [0.2in]
\epc \mbox{and} \hspace{0.7in} v_0, \ v_s \, \geq \, 0, \epc s \, = \, 1, \cdots, S \\ [5pt]
= \, \mbox{CVaR}_{\alpha}( f(x,\wt{\omega}) ) + \displaystyle{
\operatornamewithlimits{\mbox{maximum}}_v
} \epc \displaystyle{
\sum_{s=1}^S
} \, \left[ \, z_s - \mbox{CVaR}_{\alpha}( f(x,\wt{\omega}) ) \right] \, v_s \\ [0.2in]
\hspace{1.3in} \mbox{subject to} \epc \displaystyle{
\frac{-p_s}{1 - \alpha}
} \, \left[ \, \displaystyle{
\sum_{s^{\, \prime} = 1}^S
} \, v_{s^{\, \prime}} - 1 \, \right] + v_s \, \leq \, 0, \epc s \, = \, 1, \cdots, S \\ [0.2in]
\hspace{1.3in} \mbox{and} \hspace{0.5in} v_s \, \geq \, 0, \epc s \, = \, 1, \cdots, S,
\end{array} \]
where the last equality follows by the substitution: $v_0 = \displaystyle{
\sum_{s=1}^S
} \, v_s - 1$ and by noticing that the nonnegativity of $v_0$ can be dropped, provided that $\{ v_s \}_{s=1}^S$ belongs
to the fixed set:
\[
{\cal W} \, \triangleq \, \left\{ \, v \, \in \, \mathbf{R}^S_+ \, \mid \, \underbrace{\displaystyle{
\frac{-p_s}{1 - \alpha}
} \, \left[ \, \displaystyle{
\sum_{s^{\, \prime} = 1}^S
} \, v_{s^{\, \prime}} - 1 \, \right] + v_s \, \leq \, 0, \epc s \, = \, 1, \cdots, S}_{\mbox{can be written as
$\left[ \, \Id_S - \displaystyle{
\frac{1}{1 - \alpha}
} \, p \onebld_S^T \, \right] v + \displaystyle{
\frac{p}{1 - \alpha}
} \, \leq \, 0$}} \right\},
\]
where $\Id_S$ is the identity matrix of order $S$, $p$ is the $S$-vector of probabilities $\{ p_s \}_{s=1}^S$,
and $\onebld_S$ is the $S$-vector of ones.  Let $\left\{ v^j \triangleq ( v^j_s )_{s=1}^S \right\}_{j=1}^J$ be
the finite family of extreme points of ${\cal W}$ for some integer $J > 0$.  We then obtain the expression:
\begin{equation} \label{eq:our VaR formula} 
\mbox{VaR}_{\alpha}( f(x,\wt{\omega}) ) \, = \, \mbox{CVaR}_{\alpha}( f(x,\wt{\omega}) ) + \displaystyle{
\operatornamewithlimits{\mbox{maximum}}_{1 \leq j \leq J}
} \ \displaystyle{
\sum_{s=1}^S
} \, \left[ \, f(x,\omega^s) - \mbox{CVaR}_{\alpha}( f(x,\wt{\omega}) ) \, \right] \, v_s^j,
\end{equation}
which shows that $\mbox{VaR}_{\alpha}( f(x,\wt{\omega}) )$ is the pointwise maximum of a finite family of dc
functions.  As such it is itself a dc function of $x$, by known properties of difference-convexity.
The above expression is different from the one in \cite[page~866]{Wozabal12} that has the following form
for the case $f(x,\wt{\omega}) = x^T \wt{\omega}$:
\begin{equation} \label{eq:wozabal12}
\mbox{VaR}_{\alpha}( x^T\wt{\omega} ) \, = \, \displaystyle{
\frac{\alpha}{\gamma}
} \, \mbox{CVaR}_{\alpha}( x^T\wt{\omega} ) + \left( \, 1 -  \displaystyle{
\frac{\alpha}{\gamma}
} \, \right) \, \mbox{CVaR}_{\alpha - \gamma}( x^T\wt{\omega} )
\end{equation}
where $\gamma \in ( 0,\alpha )$ is a constant that is independent of $x$.  [The bilinearity of the portfolio return
$x^T\wt{\omega}$ is not essential.]   When the scenarios have equal probabilities, the constant $\gamma$ can easily
be determined.  In the general case of un-equal scenario probabilities, the constant $\gamma$ can still be obtained
by solving a bin packing problem.  In the considered examples in the cited reference, it was observed that a small $\gamma$
relative to the scenario probabilities was sufficient to validate the formula.
In contrast, our dc decomposition of
$\mbox{VaR}_{\alpha}( f(x,\wt{\omega}) )$ replaces the bin-packing step by the enumeration of the extreme points
of the special polyhedron ${\cal W}$.
Detailed investigation of the connection of these two formulae of VaR in terms of CVaR
and how the alternative expression (\ref{eq:our VaR formula}) can be used for algorithmic design are beyond the scope of this paper.
In the next subsection, we extend the above derivation to a piecewise linear utility based OCE.

\gap

The family $\left\{ v^j \triangleq ( v^j_s )_{s=1}^S \right\}_{j=1}^J$
of extreme points of the special set ${\cal W}$ depends only on the probabilities $\{ p_s \}_{s=1}^S$ and not on
the realizations $\{ \omega^s \}_{s=1}^S$.  Properties of these extreme points are not known at the present time;
understanding the polytope ${\cal W}$ and its extreme points could benefit the design
of efficient descent algorithms for optimizing $\mbox{VaR}_{\alpha}( f(x,\wt{\omega}) )$.  This is a worthwhile investigation
that is left for future research.  In the special case where the scenario probabilities $\{ p_s \}_{s=1}^S$ are all equal,
we expect that the formula (\ref{eq:our VaR formula}) can be significantly simplified.

\subsection{Extension to an OCE with piecewise linear utility} \label{subsec:OCE}

Consider the function $m_u( f(x,\wt{\omega}) )$ with a concave piecewise linear utility
function $u(t) \triangleq \displaystyle{\min_{1 \leq i \leq I}} \, a_i t + \alpha_i$ for some positive integer $I$ and
scalars $\{ a_i,\alpha_i \}_{i=1}^I$ with each $a_i \geq 0$.  Omitting the derivation which is similar to the above for the VaR, we can show that
\[
m_u( f(x,\wt{\omega}) ) = {\cal O}_u( f(x,\wt{\omega}) ) + \displaystyle{
\operatornamewithlimits{\mbox{minimum}}_{1 \leq j \leq J}
} \, \displaystyle{
\sum_{i=1}^I
} \, \displaystyle{
\sum_{s=1}^S
} \, \left\{ \, a_i \, \left[ \, f(x,\omega^s) - {\cal O}_u( f(x,\wt{\omega})) \, \right] + \alpha_i \, \right\} \, \varphi_{is}^j
\]
where each $\varphi^j = ( \varphi^j_{is} )_{(i,s) = 1}^{(I,S)}$ is an extreme point of a corresponding set
\[
\Phi \, \triangleq \, \left\{ \, \varphi \, \in \mathbb{R}_+^{IS} \, \mid \, \displaystyle{
\sum_{i=1}^I
} \, \displaystyle{
\sum_{s^{\prime}=1}^S
} \, \left( \, p_s \, a_i - \delta_{s^{\prime} s} \, \right) \, \varphi_{is^{\prime}} \, = \, p_s, \ \forall \, s \, = \, 1, \cdots, S \, \right\},
\]
where $\delta_{s^{\prime} s} \triangleq \left\{ \begin{array}{ll}
1 & \mbox{if $s^{\prime} = s$} \\
0 & \mbox{if $s^{\prime} \neq s$}
\end{array} \right.$.
Similar to Proposition~\ref{pr:CVaR}, the following result shows that ${\cal O}_u( f(x,\wt{\omega}))$ is a dc function of
$x$ provided that $f(\bullet,\omega)$ is.  Thus so is $m_u( f(\bullet,\wt{\omega}) )$.

\begin{proposition} \label{pr:OCE} \rm
For every $\omega \in \Omega$, let $f(x,\omega) = p(x,\omega) - q(x,\omega)$ be the dc decomposition of $f(\bullet,\omega)$
on a convex set ${\cal D} \subseteq \mathbb{R}^n$.  For a concave piecewise linear utility
function $u(t) \triangleq \displaystyle{\min_{1 \leq i \leq I}} \, a_i t + \alpha_i$ for some positive integer $I$ and
scalars $\{ a_i,\alpha_i \}_{i=1}^I$ with each $a_i \geq 0$, it holds that
\[
{\cal O}_u(f(x,\wt{\omega})) \, = \, \underbrace{\left( \, \displaystyle{
\sum_{i=1}^I
} \, a_i \, \right) \, \E p(x,\wt{\omega})}_{\mbox{cvx in $x$}} + \underbrace{\displaystyle{
\max_{\eta \in \mathbb{R}}
} \, \left[ \, \eta - \E\displaystyle{
\max_{1 \leq i \leq I}
} \, \left\{ \, \left( \, \displaystyle{
\sum_{i^{\, \prime} \neq i}
} \, a_{i^{\, \prime}} \, \right) p(x,\wt{\omega}) + a_i \left( \, q(x,\wt{\omega}) + \eta \, \right) - \alpha_i \, \right\} \, \right]}_{\mbox{cve in $x$}}.
\]
\end{proposition}

{\bf Proof.}  We have, for each $\omega \in \Omega$,
\[ \begin{array}{lll}
u( f(x,\omega) - \eta ) &  = & \displaystyle{\min_{1 \leq i \leq I}} \,
\left\{ \, a_i \, \left[ \, p(x,\omega) - q(x,\omega) - \eta \, \right] + \alpha_i \, \right\} \\ [5pt]
& = & \left( \, \displaystyle{
\sum_{i=1}^I
} \, a_i \, \right) \, p(x,\omega) - \underbrace{\displaystyle{
\max_{1 \leq i \leq I}
} \, \left\{ \, \left( \, \displaystyle{
\sum_{i^{\, \prime} \neq i}
} \, a_{i^{\, \prime}} \, \right) \, p(x,\omega) + a_i \, \left( \, q(x,\omega) + \eta \, \right) - \alpha_i \, \right\}}_{\mbox{jointly cvx in $(x,\eta)$}}.
\end{array} \]
Thus the desired expression of ${\cal O}_u(f(x,\wt{\omega}))$ follows readily.  \hfill $\Box$

\subsection{Variance, standard deviation, and exponential densities} \label{subsec:statistical functions}

We next show that the composite variance and standard deviations of dc functionals are dc, and so are the composite logarithmic
and exponential functions.  The proof of the former functions is
based on several elementary facts of dc functions, which we summarize in Lemma~\ref{lm:general dc properties} below.
While these facts all pertain to composite functions and are generally known
in the dc literature (see e.g.\ \cite[Theorem~II, page~708]{Hartman59} and \cite{BBorwein11,VeselyZajicek09}),
we give their proofs in order to highlight the respective dc decompositions of the functions in question.  Such explicit decompositions
are expected to be useful in applications.

\begin{lemma} \label{lm:general dc properties} \rm
Let ${\cal D}$ be a convex subset in $\mathbb{R}^n$.  The following statements are valid:
\begin{description}
\item[(a)] The square of a dc function on ${\cal D}$ is dc on ${\cal D}$.
\item[(b)] The product of two dc functions on ${\cal D}$ is dc on ${\cal D}$.
\item[(c)]
Let $f_i(x) =  g_i(x) - h_i(x)$ be a dc decomposition of $f_i$ on ${\cal D}$, for $i = 1, \cdots, I$.
Then $\| F(x) \|_2$ is a dc function on ${\cal D}$, where $F(x) \triangleq \left( f_i(x) \right)_{i=1}^I$.
\end{description}
\end{lemma}

\noindent {\bf Proof.} (a)  Let $f = g - h$ be a dc decomposition of $f$ with $g$ and $h$ being both nonnegative convex functions.
We have
\begin{equation} \label{eq:squared dc}
f^2 \, = \, 2 \, ( \, g^2 + h^2 \, ) - ( \, g + h \, )^2.
\end{equation}
This gives a dc representation of $f^2$ because the square of a nonnegative convex function is convex.

\gap

(b)  Let $\wh{f}$ and $\wt{f}$ be any two dc functions on ${\cal D}$, we can write
\[
\wh{f}(x) \, \wt{f}(x) \, = \, \thalf \, \left[ ( \, \wh{f}(x) + \wt{f}(x) \, )^2 - \wh{f}(x)^2 - \wt{f}(x)^2 \, \right],
\]
which shows that the product $\wh{f} \wt{f}$ is dc.

\gap

(c) Using the fact that $\| v \|_2 = \displaystyle{
\max_{u \, : \, \| u \|_2 = 1}
} \, u^Tv$, we can express
\[ \begin{array}{l}
\| \, F(x) \, \|_2 \, = \, \displaystyle{
\max_{u \, : \, \| u \|_2 = 1}
} \, \displaystyle{
\sum_{i=1}^I
} \, u_i \, ( \, g_i(x) - h_i(x) \, ) \\ [0.2in]
\epc = \, \displaystyle{
\sum_{i=1}^I
} \, \left[ \, \underbrace{-g_i(x) - h_i(x)}_{\mbox{cve function}} \, \right] +
\displaystyle{
\max_{u \, : \, \| u \|_2 = 1}
} \, \displaystyle{
\sum_{i=1}^I
} \, \left[ \, \underbrace{( \, u_i + 1 \, ) \, g_i(x) + ( \, 1 - u_i \, ) \, h_i(x)}_{\mbox{$\triangleq \varphi_i(x,u_i)$}} \, \right].
\end{array} \]
Since each $\varphi_i(\bullet,u_i)$ is a convex function and the pointwise maximum of a family of convex functions is convex, the
above identity readily gives a dc decomposition of $\| F(x) \|_2$.  \hfill $\Box$

\gap

{\bf Remark.}  It follows from the above proof in part (c) that the function $\| F(x) \|_2 + \displaystyle{
\sum_{i=1}^I
} \, \left( g_i(x) + h_i(x) \right)$ is convex.  This interesting side-result is based on a pointwise maximum representation of
the convex Euclidean norm function; the result can be extended to the composition of a convex with a dc
function by using the conjugacy
theory of convex functions; see Proposition~\ref{pr:composite log} for an illustrative result of this kind.
\hfill $\Box$

\gap

Utilizing the expression (\ref{eq:squared dc}), we give a dc decomposition of $\sigma^2(f(x,\wt{\omega}))$.

\begin{proposition} \label{pr:variance} \rm
For every $\omega \in \Omega$, let $f(x,\omega) = p(x,\omega) - q(x,\omega)$ be a dc decomposition of $f(\bullet,\omega)$
on  ${\cal D} \subseteq \mathbb{R}^n$ with $p(\bullet,\omega)$ and $q(\bullet,\omega)$ both nonnegative and convex.
Then
\begin{equation} \label{eq:variance dc}
\begin{array}{l}
\sigma^2(f(x,\wt{\omega})) \, = \, 2 \, \underbrace{\E\left[ \, p(x,\wt{\omega})^2 + q(x,\wt{\omega})^2 \, \right] +
\left\{ \, \E p(x,\wt{\omega}) + \E q(x,\wt{\omega}) \, \right\}^2}_{\mbox{cvx in $x$}} \\ [0.2in]
\hspace{1in} - \left\{ \,
\underbrace{\E\left[ \, p(x,\wt{\omega}) + q(x,\wt{\omega}) \, \right]^2 + 2
\left[ \, \left( \, \E p(x,\wt{\omega}) \, \right)^2 + \left( \, \E q(x,\wt{\omega}) \, \right)^2 \, \right]}_{\mbox{cvx in $x$}}
\, \right\} \, .
\end{array}
\end{equation}
\end{proposition}

\noindent {\bf Proof.}  We have
\[ \begin{array}{l}
\sigma^2(f(x,\wt{\omega})) \, = \, \E f(x,\wt{\omega})^2 - \left( \, \E f(x,\wt{\omega}) \, \right)^2 \\ [5pt]
= \, \E\left[ 2 \left\{ \, p(x,\wt{\omega})^2 + q(x,\wt{\omega})^2 \, \right\} -
\left\{ \, p(x,\wt{\omega}) + q(x,\wt{\omega}) \, \right\}^2 \, \right] \\ [5pt]
\hspace{1in} - \, \left[ \, 2
\left\{ \, \left( \, \E p(x,\wt{\omega}) \, \right)^2 + \left( \, \E q(x,\wt{\omega}) \, \right)^2 \, \right\}
- \left\{ \, \E p(x,\wt{\omega}) + \E q(x,\wt{\omega}) \, \right\}^2 \right],
\end{array} \]
from which (\ref{eq:variance dc}) follows readily.  \hfill $\Box$

\gap

The proof that $\sigma(f(x,\wt{\omega}))$ is dc requires $\wt{\omega}$ to be
discretely distributed and is based on the expression:
\[
\sigma(f(x,\wt{\omega})) \, = \, \sqrt{ \displaystyle{
\sum_{s=1}^S
} \, p_s \, \left[ \, \underbrace{f(x,\omega^s) - \displaystyle{
\sum_{s^{\, \prime}=1}^S
} \, p_{s^{\, \prime}} \, f(x,\omega^{s^{\, \prime}})}_{\mbox{dc in $x$}} \, \right]^2};
\]
thus, $\sigma(f(x,\wt{\omega}))$ is the 2-norm of the vector dc function
$x \mapsto \left( \sqrt{p_s} \left[ f(x,\omega^s) - \displaystyle{
\sum_{s^{\, \prime}=1}^S
} \, p_{s^{\, \prime}} \, f(x,\omega^{s^{\, \prime}}) \right] \right)_{s=1}^S$ (i.e., all its components are dc functions).
As such, the dc-property of $\mbox{std}(f(\bullet,\wt{\omega}))$ follows from part (c) of Lemma~\ref{lm:general dc properties}.

\gap

We end this section by addressing two key functions in the one-parameter exponential density estimation discussed in
Subsection~\ref{subsec:one-parameter exponential}.
In generic notation, the first
function is the composite $b \circ m(x)$ for $x$ in a convex set ${\cal D} \subseteq \mathbb{R}^n$.  The univariate function
$m$ is convex and non-decreasing, and admits a special dc decomposition
$m(x) = p(x) - \displaystyle{
\max_{1 \leq i \leq I}
} \, \left[ \, ( a^i )^Tx + \alpha_i \, \right]$, where $p$ is a convex function on ${\cal D}$ and each pair $( a^i,\alpha_i ) \in \mathbb{R}^{n+1}$;
such a function $m$ includes as a special case a piecewise affine function where $p$ is also the pointwise maximum of finitely many affine
functions.  The second function is $\log f(x)$, where $f$ is a dc function bounded away from zero.  Extensions of these deterministic functions
to the expected-value function $\E\left[ \log(f(x,\wt{\omega})) \right]$ is straightforward when the random variable is discretely distributed
so that
\[
\E\left[ \log(f(x,\wt{\omega})) \right] \, = \, \displaystyle{
\sum_{s=1}^S
} \, p_s \, \log( f(x,\omega^s) ).
\]

\begin{proposition} \label{pr:composite incr cvx} \rm
Let $b : \mathbb{R} \to \mathbb{R}$ be a convex non-decreasing function and $m(x) = p(x) - \displaystyle{
\max_{1 \leq i \leq I}
} \, \left[ \, ( a^i )^Tx + \alpha_i \, \right]$ with $p$ being a convex function on the convex set ${\cal D} \subseteq \mathbb{R}^n$.
It holds that the composite function $b \circ m$ is dc on ${\cal D}$.
\end{proposition}

{\bf Proof.}  By the non-decreasing property of $b$, we can write
\[
b \circ m(x) \, = \, \displaystyle{
\min_{1 \leq i \leq I}
} \, b\left( \, p(x) - ( a^i )^Tx - b_i \, \right),
\]
with each function $x \mapsto b\left( \, p(x) - ( a^i )^Tx - b_i \, \right)$ being convex.  The above expression shows that
the composite $b \circ m(x)$ is the pointwise minimum of finitely many convex functions; hence $b \circ m$ is dc.  \hfill $\Box$

\gap

{\bf Remark.}  While it is known that the composition of two dc functions is dc if their respective domains have some openness/closedness
properties; see e.g.\ \cite[Theorem~II]{Hartman59}.  The dc decomposition of the composite function is rather complex and not as simple as
the one in the special case of Proposition~\ref{pr:composite incr cvx}.  Whether the latter decomposition can be extended to the general
case where the pointwise (finite) maximum term in the function $m(x)$ is replaced by an arbitrary convex function is not known at the present time.
\hfill $\Box$

\gap

\begin{proposition} \label{pr:composite log} \rm
Let $\theta(x) = -\log( f(x) )$ where $f(x) = p(x) - q(x)$ is a dc function on a convex set ${\cal D} \subseteq \mathbb{R}^n$
where $\displaystyle{
\inf_{x \in {\cal D}}
} \, f(x) > 0$.  Then there exists a scalar
$M > 0$ such that  the function $\theta(x) + M \, p(x)$ is convex on ${\cal D}$;
thus $\theta$ has the dc decomposition
$\theta(x) \, = \, \left[ \, \theta(x) + M \, p(x) \, \right] - M \, p(x)$
on ${\cal D}$.
\end{proposition}

\noindent {\bf Proof.}  Let $1/M \triangleq \displaystyle{
\inf_{x \in {\cal D}}
} \, f(x)$.  Consider the conjugate function of the univariate convex function $\zeta(t) \triangleq -\log(t)$ for $t > 0$.
We have
\[
\zeta^*(v) \, \triangleq \displaystyle{
\sup_{t > 0}
} \, \left\{ \, v \, t + \log(t) \, \right\} \, = \, -1 - \log(-v), \epc \mbox{for $v \, < \, 0$}.
\]
Hence, by double conjugacy,
\[
-\log(t) \, = \, \displaystyle{
\sup_{v < 0}
} \, \left\{ \, v \, t + 1 + \log(-v) \, \right\}, \epc \mbox{for $t \, > \, 0$},
\]
where the sup is attained at $v = -1/t$.  Since $p(x) - q(x) \geq 1/M$ on ${\cal D}$, it follows that
\[ \begin{array}{lll}
-\log\left( p(x) - q(x) \right) & = & \displaystyle{
\sup_{v < 0}
} \, \left\{ \, v \, ( \, p(x) - q(x) \, ) + 1 + \log(-v) \, \right\} \\ [5pt]
& = & \displaystyle{
\sup_{-M \leq v < 0}
} \, \left\{ \, v \, ( \, p(x) - q(x) \, ) + 1 + \log(-v) \, \right\} \\ [5pt]
& = & \displaystyle{
\sup_{-M \leq v < 0}
} \, \left\{ \, \underbrace{( \, v + M \, ) \, p(x) + ( \, -v \, ) \, q(x) + 1 + \log(-v)}_{\mbox{cvx in $x$}} \, \right\} - M \, p(x).
\end{array} \]
This shows that $\theta(x) + M \, p(x)$ is equal to the pointwise maximum of a family of convex functions, hence is convex.  \hfill $\Box$


%
%

\section{A Study of the Value Function $\mbox{qp}_{\rm opt}(q,b)$} \label{sec:quadratic value fnc}

Properties of the QP-value function are clearly important for
the understanding of the quadratic recourse function $\psi(x,\omega)$, which is equal to $\mbox{qp}_{\rm opt}(q(x,\omega),b(x,\omega))$
where $q(x,\omega) \triangleq f(\omega) + G(\omega)x$ and $b(x,\omega) \triangleq \xi(\omega) - C(\omega)x$.  In the
analysis of the general QP (\ref{eq:copositive qp}), we do not make a semi-definiteness assumption on $Q$.  Instead the copositivity
of $Q$ on $D_{\infty}$ is essential because if there exists a recession vector $v$ in this
cone such that $v^TQv < 0$, then for all $b$ for which ${\cal P}_D(b) \neq \emptyset$, we have $\mbox{qp}_{\rm opt}(q,b) = -\infty$ for all
$q \in \mathbb{R}^n$.

\gap

We divide the derivation of the dc property (and an associated dc representation) of $\mbox{qp}_{\rm opt}(q,b)$ into two cases: (a) when $Q$ is
positive definite, and (b) when $Q$ is copositive on $D_{\infty}$.  In the former case, the dc decomposition
is easy to derive and much simpler, whereas the latter case is much more involved.  Indeed, when $Q$ is positive definite, we have
\begin{equation} \label{eq:pd QP}
\begin{array}{lllll}
\mbox{qp}_{\rm opt}(q,b) & = & -\thalf \, q^T Q^{-1} q + & \displaystyle{
\operatornamewithlimits{\mbox{minimum}}_{z}
} & \thalf \, \left( \, z + Q^{-1}q \, \right)^T Q \left( \, z + Q^{-1}q \, \right) \\ [7pt]
& & & \mbox{subject to} & Dz \, \geq \, b \\ [0.1in]
& = & -\thalf \, q^T Q^{-1} q + & \displaystyle{
\operatornamewithlimits{\mbox{minimum}}_y
} & \thalf \, y^T Q y \\ [7pt]
& & & \mbox{subject to} & Dy \, \geq \, b^{\, \prime} \, \triangleq \, b + DQ^{-1}q.
\end{array} \end{equation}
The lemma below shows that the minimum of the second summand in the above expression is a convex function of the right-hand vector $b^{\, \prime}$.
This yields a rather simple dc decomposition of $\mbox{qp}_{\rm opt}(q,b)$.  The proof of the lemma is easy and omitted.

\begin{lemma} \label{lm:convexity rhs} \rm
Let $f : \mathbb{R}^m \to \mathbb{R}$ be a strongly convex function and let $D \in \mathbb{R}^{k \times m}$.
The value function:
\[ \begin{array}{lll}
\varphi(b) \, \triangleq & \displaystyle{
\operatornamewithlimits{\mbox{minimum}}_y
} & f(y) \\ [3pt]
& \mbox{subject to} & Dy \, \geq \, b
\end{array} \]
is a convex function on the domain $\mbox{Range}(D) - \mathbb{R}^k_+$ which consists of all vectors $b$ for which there exists $y \in \mathbb{R}^m$
satisfying $Dy \geq b$.  \hfill $\Box$
\end{lemma}

It is worthwhile to point out that while it is fairly easy to derive the above ``explicit'' dc decomposition of $\mbox{qp}_{\rm opt}(q,b)$ when $Q$ is
positive definite, the same cannot be said when $Q$ is positive semidefinite.
Before proceeding further, we refer the reader to the monograph \cite{LTYen05} for an extensive study of
indefinite quadratic programs.  Yet results therein do not provide clear descriptions of the structure of
(a) the optimal solution set of the QP (\ref{eq:copositive qp}) for fixed $(q,b)$ and (b) the domain $\mbox{dom}(Q,D)$ of finiteness,
and (c) the optimal value function $\mbox{qp}_{\rm opt}(q,b)$, all when $Q$ is an indefinite matrix.  Instead
we rely on an early result \cite{GiannessiTomasin73} and the recent study \cite{HuMitchellPang12} to formally state and prove
these desired properties
of the QP (\ref{eq:copositive qp}) when $Q$ is copositive on the recession cone of the feasible region.  This is the main content of
Proposition~\ref{pr:copositive results} below.  The key to this proposition
is the following consequence of a well-known property of a quadratic program which was originally due to Frank-Wolfe and subsequently
refined by Eaves \cite{Eaves71}; see \cite[Theorem~2.2]{LTYen05}.  In the lemma below and also subsequently, the
$\perp$ notation denotes perpendicularity, which expresses the complementarity condition between two nonnegative vectors.

\begin{lemma} \label{lm:Eaves} \rm
Suppose that $Q$ is copositive on $D_{\infty}$.  A pair $(q,b) \in \mbox{dom}(Q,D)$ if and only if ${\cal P}_D(b) \neq \emptyset$ and
\begin{equation} \label{eq:copositive boundedness}
\left. \begin{array}{r}
Qv - D^{\, T}\eta \, = \, 0 \\ [5pt]
0 \, \leq \, \eta \, \perp \, Dv \, \geq \, 0 \\ [5pt]
Dz \, \geq \, b
\end{array} \right\} \ \Rightarrow \ ( \, q + Qz \, )^Tv \, \geq \, 0.
\end{equation}
\end{lemma}

\noindent {\bf Proof.}  By \cite[Theorem~2.2]{LTYen05}, $(q,b) \in \mbox{dom}(Q,D)$ if and only if ${\cal P}_D(b) \neq \emptyset$ and
\begin{equation} \label{eq:copositivity=0}
\left. \begin{array}{r}
v^TQv \, = \, 0, \epc Dv \, \geq \, 0 \\ [5pt]
Dz \, \geq \, b
\end{array} \right\} \ \Rightarrow \ ( \, q + Qz \, )^Tv \, \geq \, 0.
\end{equation}
By the copositivity of $Q$ on $D_{\infty}$, we have
\[ \begin{array}{lll}
\left[ \, v^TQv \, = \, 0, \epc Dv \, \geq  \, 0 \, \right] & \Leftrightarrow &
v \mbox{ is a minimizer of } \thalf \, u^TQu \mbox{ for } u \, \in \, D_{\infty} \\ [5pt]
& \Leftrightarrow &
\exists \, \eta \mbox{ such that }
\left\{ \begin{array}{r}
Qv - D^{\, T}\eta \, = \, 0 \\ [5pt]
0 \, \leq \, \eta \, \perp \, Dv \, \geq \, 0.
\end{array} \right.
\end{array} \]
Thus (\ref{eq:copositivity=0}) and (\ref{eq:copositive boundedness}) are equivalent.  \hfill $\Box$

\gap

Interestingly, while the result below is well known in the case of a positive semidefinite $Q$, its extension
to a copositive $Q$ turns out to be not straightforward, especially under no boundedness assumption whatsoever.
There is an informal assertion, without proof or citation, of the result below (in particular, part (c))
for a general qp without the copositive assumption \cite[page~88]{LuoPangRalph96}.  This is a careless oversight
on the authors' part as the proof below is actually quite involved.

\begin{proposition} \label{pr:copositive results} \rm
Suppose that $Q$ is copositive on $D_{\infty}$.  The following statements hold:
\begin{description}
\item[\rm (a)] $\mbox{qp}_{\rm opt}(q,b)$ is continuous on $\mbox{dom}(Q,D)$;
\item[\rm (b)] $\mbox{dom}(Q,D)$ is the union of finitely many polyhedra (thus is a closed set);
\item[\rm (c)] there exist a finite family ${\cal F} \triangleq \{ S_F \}$ of polyhedra in $\mathbb{R}^{m+k}$ and finitely many quadratic functions $\{ \mbox{qp}_F \}$ such that $\mbox{qp}_{\rm opt}(q,b) =
\displaystyle{
	\min_{F : (q,b) \in S_F}
} \mbox{qp}_F(q,b)$; hence $\mbox{qp}_{\rm opt}(q,b)$ is
a piecewise quadratic function on $\mbox{dom}(Q,D)$;
\item[\rm (d)] for each pair $(q,b) \in \mbox{dom}(Q,D)$, the optimal solution set of the QP (\ref{eq:copositive qp}) is the union of finitely
many polyhedra.
\end{description}
\end{proposition}

\noindent {\bf Proof.}
(a) Let $\{ (q^{\nu},b^{\nu}) \} \subseteq \mbox{dom}(Q,D)$ be a sequence of vectors converging to the pair $( q^{\infty},b^{\infty} )$.
For each $\nu$, let $z^{\nu} \in  \mbox{qp}_{\rm sol}(q^{\nu},b^{\nu})$.  Note that this sequence $\{ z^{\nu} \}$ need not be bounded.  Nevertheless,
by the renowned Hoffman error bound for linear inequalities \cite[Lemma~3.2.3]{FacchineiPang03}, it follows that the following two properties
hold: (a) ${\cal P}_D(b^{\infty}) \neq \emptyset$ and
there exists a sequence $\{ \bar{z}^{\, \nu} \} \subseteq {\cal P}_D(b^{\infty})$ such that $\{ \bar{z}^{\nu} - z^{\nu} \} \to 0$; and (b) for every
$z \in {\cal P}_D(b^{\infty})$, there exists, for every $\nu$ sufficiently large, a vector $y^{\nu} \in {\cal P}_D(b^{\nu})$ such that the sequence
$\{ y^{\nu} \}$ converges to $z$.   Hence, for every $( \eta,v )$ satisfying the left-hand conditions in (\ref{eq:copositive boundedness}),
we have
\[
( \, q^{\infty} + Qz \, )^Tv \, = \, \displaystyle{
\lim_{\nu \to \infty}
} \, ( \, q^{\nu} + Qy^{\nu} \, )^Tv \, \geq \, 0.
\]
Hence $\mbox{qp}_{\rm sol}(q^{\infty},b^{\infty}) \neq \emptyset$; thus
$( q^{\infty},b^{\infty} ) \in \mbox{dom}(Q,D)$.  Moreover, by taking
the vector $z \in \mbox{qp}_{\rm sol}(q^{\infty},b^{\infty})$,
we have
\[
\mbox{qp}_{\rm opt}(q^{\infty},b^{\infty}) \, = \,
\zeta(z) \, = \, \displaystyle{
\lim_{\nu \to \infty}
} \, \zeta(y^{\nu}) \, \geq \, \displaystyle{
\limsup_{\nu \to \infty}
} \, \zeta(z^{\nu}) \, = \, \displaystyle{
\limsup_{\nu \to \infty}
} \, \zeta( \bar{z}^{\nu}) \, \geq \, \displaystyle{
\liminf_{\nu \to \infty}
} \, \zeta( \bar{z}^{\nu}) \, \geq \, \mbox{qp}_{\rm opt}(q^{\infty},b^{\infty})
\]
from which it follows that $\displaystyle{
\lim_{\nu \to \infty}
} \, \mbox{qp}_{\rm opt}(q^{\nu},b^{\nu}) = \displaystyle{
\lim_{\nu \to \infty}
} \, \zeta(z^{\nu})$ exists and equals to $\mbox{qp}_{\rm opt}(q^{\infty},b^{\infty})$.  Thus (a) holds.

\gap

To prove (b), we use an early result of quadratic programming \cite{GiannessiTomasin73} stating that for a pair
$(q,b)$ in $\mbox{dom}(Q,D)$, the value $\mbox{qp}_{\rm opt}(q,b)$ is equal to the minimum of the quadratic objective function
$\zeta(z)$ on the set of stationary solutions of the problem, and also use the fact \cite{LuoTseng92} that the set of values
of $\zeta$ on the set of stationary solutions is finite.  More explicitly,  the stationarity conditions of (\ref{eq:copositive qp}),
or Karush-Kuhn-Tucker (KKT) conditions, are given by the following mixed complementarity conditions:
\begin{equation} \label{eq:KKT QP}
\begin{array}{r}
q + Qz - D^{\, T}\eta \, = \, 0 \\ [5pt]
0 \, \leq \, \eta \, \perp \, Dz - b \, \geq \, 0.
\end{array} \end{equation}
In turn the above conditions can be decomposed
into a finite, but exponential, number of linear inequality systems by considering the index subsets
${\cal I} \subseteq \{ 1, \cdots, k \}$ each with complement ${\cal J}$ such that
\begin{equation} \label{eq:KKT-calI} \begin{array}{lll}
q + Qz - ( D_{{\cal I} \bullet} )^T \eta_{\cal I} & = & 0, \\ [5pt]
D_{{\cal I} \bullet}z - b_{\cal I} \, = \, 0 & \leq & \eta_{\cal I} \\ [5pt]
D_{{\cal J} \bullet}z - b_{\cal J} \, \geq \, 0 & = & \eta_{\cal J}.
\end{array} \end{equation}
The above linear inequality system defines a polyhedral set in $\mathbb{R}^{m+k}$, which we denote $\mbox{KKT}({\cal I})$.
Alternatively, a tuple $( q,b )$ satisfies (\ref{eq:KKT-calI}) if and only if 
\[
\left( \begin{array}{l}
q \\ [5pt]
b_{\cal I} \\ [5pt]
b_{\cal J}
\end{array} \right) \, \in \, \left[ \begin{array}{ccc}
-Q & ( D_{{\cal I} \bullet} )^T & 0\\ [5pt]
D_{{\cal I} \bullet} & 0 & 0 \\ [5pt]
D_{{\cal J} \bullet} & 0 & -\Id_{| {\cal J} |}
\end{array} \right] \left( \begin{array}{c}
\mathbb{R}^m \\ [5pt]
\mathbb{R}_+^{| {\cal I} |} \\ [5pt]
\mathbb{R}_+^{| {\cal J} |}
\end{array} \right) \, \triangleq \, \underbrace{\mbox{QD}({\cal I})}_{\mbox{a polyhedral cone}}.
\]
The recession cone of $\mbox{KKT}({\cal I})$, denoted $\mbox{KKT}({\cal I})_{\infty}$, is a polyhedral cone
defined by the following homogeneous system in the variables $( v,\xi )$:
\[ \begin{array}{lll}
Qv - ( D_{{\cal I} \bullet} )^T \xi_{\cal I} & = & 0, \\ [5pt]
D_{{\cal I} \bullet}v \, = \, 0 & \leq & \xi_{\cal I} \\ [5pt]
D_{{\cal J} \bullet}v \, \geq \, 0 & = & \xi_{\cal J}.
\end{array} \]
This cone is the conical hull of a finite number of generators, which we denote
$\left\{ \left( v^{{\cal I},\ell},\xi^{{\cal I},\ell} \right) \right\}_{\ell=1}^{L_{\cal I}}$ for some integer
$L_{\cal I} > 0$.  Notice that these generators depend only on the pair $(Q,D)$.  In terms of them, the
implication (\ref{eq:copositive boundedness}) is equivalent to
\[
Dz \, \geq \, b \ \Rightarrow \ ( \, q + Qz \, )^Tv^{{\cal I},\ell} \, \geq \, 0, \epc \forall \, \ell \, = \, 1, \cdots, L_{\cal I} \mbox{ and all subsets
${\cal I}$ of $\{ 1, \cdots, k \}$}.
\]
In turn, for a fixed pair $( {\cal I},\ell )$, the latter implication holds if and only if
\[ \begin{array}{llllll}
0 & \leq & q^Tv^{{\cal I},\ell} + & \displaystyle{
\operatornamewithlimits{\mbox{minimum}}_{Dz \, \geq \, b}
} & \left( \, Qv^{{\cal I},\ell} \, \right)^Tz \\ [0.1in]
& = & q^Tv^{{\cal I},\ell} + & \displaystyle{
\operatornamewithlimits{\mbox{maximize}}_{\mu \geq 0}
} & b^{\, T}\mu \hspace{1in} \mbox{by linear programming duality} \\ [5pt]
& & & \mbox{subject to} & D^T\mu \, = \, Qv^{{\cal I},\ell}.
\end{array} \]
Let ${\cal E}({\cal I},\ell) \triangleq \{ \mu^{{\cal I},\ell,\kappa} \}_{\kappa=1}^{E({\cal I},\ell)}$ be the finite set of extreme
points of the polyhedron
$\left\{ \mu \in \mathbb{R}^k_+ \mid  D^T\mu \, = \, Qv^{{\cal I},\ell} \right\}$ for some positive integer
$E({\cal I},\ell)$.  Again these extreme points depend 
on the pair $(Q,D)$ only.  It then follows that (\ref{eq:copositive boundedness})
holds if and only if for every pair $( {\cal I},\ell )$ with ${\cal I}$ being a subset of $\{ 1, \cdots, k \}$ and
$\ell = 1, \cdots, L_{\cal I}$, there exists $\mu^{{\cal I},\ell,\kappa} \in {\cal E}({\cal I},\ell)$ such that
$0 \leq q^Tv^{{\cal I},\ell} + b^{\, T} \mu^{{\cal I},\ell,\kappa}$.  Together with the feasibility condition
$b \in \mbox{range}(D) - \mathbb{R}^k_+$,
each of the latter inequalities in $(q,b)$ defines a polyhedron that defines a piece of $\mbox{dom}(Q,D)$.
Specifically let ${\cal V} \triangleq \displaystyle{
\prod_{( {\cal I},\ell )}
} \, {\cal E}({\cal I},\ell)$ where the Cartesian product ranges over all pairs $( {\cal I},\ell )$
with ${\cal I}$ being a subset of $\{ 1, \cdots, k \}$ and $\ell = 1, \cdots, L_{\cal I}$.  An element
$\boldsymbol{\mu} \in {\cal V}$ is a tuple of extreme points $\mu^{{\cal I},\ell,\kappa}$ over all pairs
$( {\cal I},\ell )$ as specified with $\kappa$ being one element in $\{ 1, \cdots, E({\cal I},\ell) \}$ corresponding
to the given $\boldsymbol{\mu}$;
every such element $\boldsymbol{\mu}$ defines a system of $\displaystyle{
\sum_{{\cal I} \subseteq \{ 1, \cdots, k \}}
} \, | \, L_{\cal I} \, |$ linear inequalities each of the form:
\[
q^Tv^{{\cal I},\ell} + b^{\, T} \mu^{{\cal I},\ell,\kappa} \, \geq \, 0,
\]
Consequently, it follows that
\[
\mbox{dom}(Q,D) \, = \, \displaystyle{
\bigcup_{\boldsymbol{\mu} \in {\cal V}}
} \, \left\{ \, ( \, q,b \, ) \, \in \, \mathbb{R}^m \, \times \, \left[ \, D\mathbb{R}^m - \mathbb{R}^k_+ \, \right] \, \mid \,
\underbrace{q^Tv^{{\cal I},\ell} + b^{\, T} \mu^{{\cal I},\ell,\kappa} \, \geq \, 0}_{\mbox{there are $\displaystyle{
\sum_{{\cal I} \subseteq \{ 1, \cdots, k \}}
} \, | \, L_{\cal I} \, |$ of these}} \, \right\},
\]
proving that $\mbox{dom}(Q,D)$ is a union of finitely (albeit potentially exponentially) many polyhedra, each denoted
by $\mbox{dom}^{\boldsymbol{\mu}}(Q,D)$ for $\boldsymbol{\mu} \in {\cal V}$.

\gap

To prove (c), let ${\cal F}$ be the family of polyhedra
\[
\left\{ \, \mbox{dom}^{\boldsymbol{\mu}}(Q,D) \, \cap \, \mbox{QD}({\cal I}^{\, \prime}) \, \mid \,
\boldsymbol{\mu} \, \in \, {\cal V} \mbox{ and } {\cal I}^{\, \prime} \, \subseteq \, \{ 1, \cdots, k \} \, \right\}.
\]
For each polyhedron in this family,
the system KKT$({\cal I}^{\, \prime})$ has a solution $( z,\eta_{{\cal I}^{\, \prime}} )$ satisfying:
\begin{equation} \label{eq:QDset}
\begin{array}{lll}
q + Qz - ( D_{{\cal I}^{\, \prime} \bullet} )^T \eta_{{\cal I}^{\, \prime}} & = & 0, \\ [5pt]
D_{{{\cal I}^{\, \prime}} \bullet}z - b_{{\cal I}^{\, \prime}} \, = \, 0 & \leq & \eta_{{\cal I}^{\, \prime}} \\ [5pt]
D_{{{\cal J}^{\, \prime}} \bullet}z - b_{{\cal J}^{\, \prime}} \, \geq \, 0, 
\end{array} \end{equation}
where ${\cal J}^{\, \prime}$ is the complement of ${\cal I}^{\, \prime}$ in $\{ 1, \cdots, k \}$.   The system (\ref{eq:QDset})
can be written as:
\[ \left[ \begin{array}{ccc}
-Q & ( D_{{\cal I}^{\, \prime} \bullet} )^T & 0 \\ [5pt]
D_{{{\cal I}^{\, \prime}} \bullet} & 0 & 0 \\ [5pt]
D_{{{\cal J}^{\, \prime}} \bullet} & 0 & -\Id_{| {\cal J}^{\, \prime} | }
\end{array} \right] \left( \begin{array}{c}
z \\ [5pt]
\eta_{{\cal I}^{\, \prime}} \\ [5pt]
s_{{\cal J}^{\, \prime}}
\end{array} \right) \, = \, \left( \begin{array}{c}
q \\ [5pt]
b_{{\cal I}^{\, \prime}} \\ [5pt]
b_{{\cal J}^{\, \prime}}
\end{array} \right), \epc ( \, \eta_{{\cal I}^{\, \prime}}, s_{{\cal J}^{\, \prime}} \, ) \, \geq \, 0.
\]
By elementary linear-algebraic operations similar to the procedure of obtaining a basic feasible solution
from an arbitrary feasible solution to a system of linear inequalities in linear programming, it follows that
there exists a matrix $M({\cal I}^{\, \prime})$,
dependent on the pair $(Q,D)$ only, such that (\ref{eq:QDset}) has a solution
with $z = M({\cal I}^{\, \prime})\left( \begin{array}{c}
q \\
b
\end{array} \right)$; i.e., (\ref{eq:QDset}) has a solution that is linearly dependent on the pair $(q,b)$.  It is not
difficult to show that the objective
function $\zeta(z)$ is a constant on the set of solutions of (\ref{eq:QDset}) [in fact, this constancy property is
the source of the finite number of values attained by the quadratic function $\zeta$ on the set of stationary solutions
of the QP (\ref{eq:copositive qp})], it follows that this constant must be a quadratic function of the pair $(q,b)$.  Since for
a fixed pair $(q,b) \in \mbox{dom}(Q,D)$,
$\mbox{qp}_{\rm opt}(q,b)$ is the minimum of these constants over all polyhedra in the family ${\cal F}$ that contains the
pair $(q,b)$, part (c) follows.

\gap

Lastly, to prove (d), it suffices to note that for a given pair $(q,b) \in \mbox{dom}(Q,D)$, the solution set of the
QP (\ref{eq:copositive qp}) is equal to set of vectors $z$ for which there exists $\eta$ such that the pair $(z,\eta)$ satisfies
the linear equality system $\mbox{KKT}({\cal I})$ and $\thalf \, ( q^Tz + b^T \eta ) = \mbox{qp}_{\rm opt}(q,b)$ for some
index subset ${\cal I}$ of $\{ 1, \cdots, k \}$.  \hfill $\Box$

\subsection{Proof the dc property}

Based on Proposition~\ref{pr:copositive results} and known properties of piecewise dc functions in general,
the following partial dc property of the value function $\mbox{qp}_{\rm opt}(q,b)$ can be proved.

\begin{corollary} \label{co:consequence of mixing} \rm
Suppose that $Q$ is copositive on $D_{\infty}$.  The value function $\mbox{qp}_{\rm opt}(q,b)$ is dc on
any open convex set contained in $\mbox{dom}(Q,D)$.
\end{corollary}

\noindent {\bf Proof.}  This follows from the mixing property of piecewise dc functions; see \cite[Proposition~2.1]{BBorwein11}
which has its source from \cite[Lemma~4.8]{VeselyZajicek89}.  \hfill $\Box$

\gap

Since $\mbox{dom}(Q,D)$ is a closed set, the above corollary does not yield the difference-convexity of $\mbox{qp}_{\rm opt}(q,b)$
on its full domain.  We give a formal statement of this property in Theorem~\ref{th:full dc} under the assumption that
$\mbox{dom}(Q,D)$ is convex.  The proof of this desired dc property is based on a more general Proposition~\ref{pr:Maher PQ}
pertaining to piecewise LC$^{\, 1}$ functions that does not
require the pieces to be quadratic functions nor the polyehdrality of the sub-domains.

\gap

Recall that a function $\theta$ is LC$^1$ (for Lipschitz continuous gradient) on an open set ${\cal O}$
in $\mathbb{R}^n$
if $\theta$ is differentiable and its gradient is a Lipschitz continuous function on ${\cal O}$.  It is
known that a LC$^1$ function is dc with the decomposition $\theta(x) = \left( \theta(x) + \displaystyle{
\frac{M}{2}
} \, x^Tx \right) - \displaystyle{
\frac{M}{2}
} \, x^Tx $, where $M$ is a positive scalar larger than the Lipschitz modulus of the gradient function $\nabla \theta$.
Another useful fact \cite[{\bf 3.2.12}]{OrtegaRheinboldt00} about LC$^{\, 1}$ functions is the inequality (\ref{eq:LC1})
that is a consequence of a mean-value function theorem of multivariate functions: namely, if $\theta$ is
LC$^{\, 1}$ with $L > 0$ being a Lipschitz modulus of $\nabla \theta$ on an open convex convex set ${\cal O}$, then
for any two vectors $x$ and $y$ in ${\cal O}$,
\begin{equation} \label{eq:LC1}
\theta(x) - \theta(y) \, \leq \, \nabla \theta(y)^T(x - y) + \displaystyle{
\frac{L}{2}
} \, \| \, x - y \, \|_2^2.
\end{equation}

\begin{proposition} \label{pr:Maher PQ} \rm
Let $\theta(x)$ be a continuous function on a convex set
${\cal S} \triangleq \displaystyle{
\bigcup_{i=1}^I
} \, S^{\, i}$ where each $S^{\, i}$ is a closed convex set in $\mathbb{R}^N$.  Suppose there exist LC$^1$
functions $\{ \theta_i(x) \}_{i=1}^I$ defined on an open set ${\cal O}$ containing ${\cal S}$
such that $\theta(x) = \theta_i(x)$ for all $x \in S^{\, i}$ and that each difference function
$\theta_{ji}(x) \triangleq \theta_j(x) - \theta_i(x)$ has dc gradients on ${\cal S}$.
It holds that $\theta$ is dc on ${\cal S}$.
\end{proposition}

\noindent {\bf Proof.}
Let $L_{ji}$ be the Lipschitz modulus of $\nabla \theta_{ji}$ on ${\cal O}$
and let $L_i \triangleq \displaystyle{
\max_{1 \leq j \leq I}
} \, L_{ji}$.  Let $\mbox{dist}(x;S^{\, i}) \triangleq \displaystyle{
\min_{z \in S^{\, i}}
} \, \| \, x - z \, \|_2 = \| \, x - \Pi_{S^{\, i}}(x) \, \|_2$ be the distance function to the set $S^{\, i}$ with $\Pi_{S^{\, i}}(x)$ being the
Euclidean projection (i.e., closest point) of the vector $x$ onto $S^{\, i}$.  We note that both $\mbox{dist}(\bullet;S^{\, i})$
and $\left[ \, \mbox{dist}(\bullet;S^{\, i}) \, \right]^2$ are convex functions.  Define
\[
\phi_i(x) \, \triangleq \, \underbrace{\mbox{dist}(x;S^{\, i})}_{\mbox{cvx in $x$}} \,
\underbrace{\displaystyle{
\max_{1 \leq j \leq I}
} \, \| \nabla \theta_{ji}(x) \|_2}_{\mbox{dc in $x$}}
+ \displaystyle{
\frac{3 \, L_i}{2}
} \, \left[ \, \mbox{dist}(x;S^{\, i}) \, \right]^2, \epc x \, \in \, \mathbb{R}^N
\]
and let $\psi_i(x) \triangleq \, \theta_i(x) + \phi_i(x)$.  The first summand,
$\mbox{dist}(x;S^{\, i}) \, \displaystyle{
\max_{1 \leq j \leq I}
} \, \| \nabla \theta_{ji}(x) \|_2$, being the product of two dc functions is dc, by Lemma~\ref{lm:general dc properties}.
Hence each function $\psi_i$ is dc.  We claim that
\begin{equation} \label{eq:dc PQ}
\theta(x) \, = \, \displaystyle{
\min_{1 \leq i \leq I}
} \, \psi_i(x), \epc \forall \, x \, \in \, {\cal S}.
\end{equation}
Since the pointwise minimum of finitely many dc functions is dc \cite[Proposition~4.1]{Tuy16}, the expression (\ref{eq:dc PQ}) is
enough to show that $\theta$ is dc on ${\cal S}$.  In turn since $\psi_i(x) = \theta_i(x)$
for all $x \in S^{\, i}$,  to show (\ref{eq:dc PQ}), it suffices to show that $\psi_i(x) \geq \theta(x)$
for all $x \in {\cal S} \setminus S^{\, i}$.  For such an $x$, let $\bar{x}^i \triangleq \Pi_{S^{\, i}}(x)$.  Since ${\cal S}$ is convex, the line
segment joining $x$ and $\bar{x}^i$ is contained in ${\cal S}$.  Hence, there exists a finite partition of the interval $[ 0,1 ]$:
\[
0 \, = \, \tau_0 \, < \, \tau_1 \, < \, \cdots \, < \, \tau_T \, < \, \tau_{T+1} \, = \, 1
\]
and corresponding indices $i_t \in \{ 1, \cdots, I \}$
for $t = 1, \cdots, T+1$ so that $\theta(x) = \theta_{i_t}(x)$ for all $x$ in the (closed) sub-segment joining
$x^{t-1} \triangleq x + \tau_{t-1} ( \bar{x}^i - x )$ to $x^t \triangleq x + \tau_t ( \bar{x}^i - x )$.  We have
\[ \begin{array}{l}
\theta(x) - \theta_i(x) \, = \, \theta(x) - \theta(\bar{x}^i) + \theta_i(\bar{x}^i) - \theta_i(x) \epc
\mbox{because $\theta(\bar{x}^i) = \theta_i(\bar{x}^i)$} \\ [0.15in]
= \, \displaystyle{
\sum_{t=1}^{T+1}
} \, \left[ \, \theta_{i_t}(x^{t-1}) - \theta_{i_t}(x^t) \, \right] - \displaystyle{
\sum_{t=1}^{T+1}
} \, \left[ \, \theta_i(x^{t-1}) - \theta_i(x^t) \, \right] \\ [0.15in]
\hspace{1in} \mbox{because $\theta(x^{t-1}) = \theta_{i_t}(x^{t-1})$ and $\theta(x^t) = \theta_{i_t}(x^t)$} \\ [0.15in]
= \, \displaystyle{
\sum_{t=1}^{T+1}
} \, \left[ \, \theta_{i_t}(x^{t-1}) - \theta_i(x^{t-1})  \, \right] - \displaystyle{
\sum_{t=1}^{T+1}
} \, \left[ \, \theta_{i_t}(x^t) - \theta_i(x^t) \, \right] \\ [0.25in]
= \, \displaystyle{
\sum_{t=1}^{T+1}
} \, \left[ \, \theta_{i_t i}(x^{t-1}) - \theta_{i_t i}(x^t) \, \right] \\ [0.25in]
\leq \, \displaystyle{
\sum_{t=1}^{T+1}
} \, \left[ \, \nabla \theta_{i_t i}(x^t)^T( \, x^{t-1} - x^t \, ) + \displaystyle{
\frac{L_{i_t i}}{2}
} \, \| \, x^{t-1} - x^t \, \|_2^2 \, \right] \epc \mbox{(by (\ref{eq:LC1}))} \\ [0.2in]
\leq \, \displaystyle{
\sum_{t=1}^{T+1}
} \, \left\{ \, \nabla \theta_{i_ti}(x)^T( \, x^{t-1} - x^t \, ) + \left[ \, \nabla \theta_{i_t i}(x^t) - \nabla \theta_{i_t i}(x) \right]^T
( \, x^{t-1} - x^t \, ) \, \right\} + \displaystyle{
\frac{\displaystyle{
\max_{1 \leq j \leq I}
} \, L_{ji}}{2}
} \, \displaystyle{
\sum_{t=1}^{T+1}
} \, \| \, x^{t-1} - x^t \, \|_2^2 \\ [0.2in]
\leq \, \displaystyle{
\max_{1 \leq j \leq I}
} \, \| \, \nabla \theta_{ji}(x) \, \|_2 \, \displaystyle{
\sum_{t=1}^{T+1}
} \, \| \, x^{t-1} - x^t \, \|_2 + \displaystyle{
\max_{1 \leq j \leq I}
} \, L_{ji} \, \| \, x - \bar{x}^i \, \|_2 \, \displaystyle{
\sum_{t=1}^{T+1}
} \, \| \, x^{t-1} - x^t \, \|_2 + \displaystyle{
\frac{L_i}{2}
} \, \| \, x - \bar{x}^i \, \|_2^2 \\ [0.2in]
= \, \displaystyle{
\max_{1 \leq j \leq I}
} \, \| \, \nabla \theta_{ji}(x) \, \|_2 \, \| \, x - \bar{x}^i \, \|_2 + \displaystyle{
\frac{3 \, L_i}{2}
} \, \| \, x - \bar{x}^i \, \|_2^2 \, = \, \phi_i(x) \epc \mbox{because $\| \, x - \bar{x}^i \, \|_2 = \mbox{dist}(x;S^{\, i})$},
\end{array} \]
where the last inequality holds by the Cauchy-Schwartz inequality, the Lipschitz property of $\nabla \theta_{i_t i}$, and the
identity $\displaystyle{
\sum_{t=1}^{T+1}
} \, \|  x^{t-1} - x^t \|_2 = \| x - \bar{x}^i \|$.
Thus $\psi_i(x) \geq \theta(x)$ for all $x \in {\cal{S}} \setminus S^{\, i}$ as claimed, and (\ref{eq:dc PQ}) follows readily.
\hfill $\Box$

\gap

We are now ready to formally state and prove the following main result of this section.

\begin{theorem} \label{th:full dc} \rm
Suppose that $Q$ is copositive on $D_{\infty}$ and that $\mbox{dom}(Q,D)$ is convex.
Then the value function $\mbox{qp}_{\rm opt}(q,b)$ is dc on $\mbox{dom}(Q,D)$, provided that $\mbox{qp}_{\rm opt}$ is a quadratic function on each polyhedral member in the family ${\cal F}$ described in part~c of Proposition~\ref{pr:copositive results}. In particular, $\mbox{qp}_{\rm opt}(q,b)$ is dc on $\mbox{dom}(Q,D)$, if $Q$ is positive semidefinite.

\end{theorem}

\noindent {\bf Proof.}  This follows readily from Proposition~\ref{pr:Maher PQ} by letting
each pair $(S_i,\theta_i)$ be the quadratic piece identified in part (c) of
Proposition~\ref{pr:copositive results}.  \hfill $\Box$

\gap

Several remarks about the above results are in order. First, we refer to two examples in the literature that show, respectively, a constraint-only \cite{Klatte85}, or an objective-only \cite{LeeTamYen05} perturbed (nonconvex) quadratic program (QP) may not be piecewise linear-quadratic; i.e. the condition of $\mbox{qp}_{\rm opt}$ being a quadratic function on each polyhedral member in the family ${\cal F}$ may not be satisfied. However, we should point out that while not satisfying the quadratic condition in the above theorem, the two numerical examples in \cite{Klatte85,LeeTamYen05} can be verified
to be dc, by the mixing property Corollary~\ref{co:consequence of mixing}.  For the one in \cite{LeeTamYen05} the domain of the value function is $\mathbb{R}^2$; for the other one \cite{Klatte85}, the domain of the function, which is restricted to the (closed) fourth quadrant in $\mathbb{R}^2$ in the reference, can be easily extended to an open convex set containing this quadrant. Thus, the two examples in \cite{Klatte85,LeeTamYen05} do not constitute a counterexample for which the value function $\mbox{qp}_{\rm opt}(q,b)$ is not dc. The question of whether we can relax the quadratic condition in Theorem~\ref{th:full dc} and still show that the value function $\mbox{qp}_{\rm opt}(q,b)$ is dc remains open. Second, when each $\theta_i(x) = x^Ta^i + \alpha_i$ is an affine function for some
$N$-vector $a^i$ and scalar $\alpha_i$, the representation (\ref{eq:dc PQ}) of the function $\theta$ becomes
\[
\theta(x) \, = \, \displaystyle{
\min_{1 \leq i \leq I}
} \, \left[ \, \theta_i(x) + \mbox{dist}(x;S^{\, i}) \, \displaystyle{
\max_{1 \leq j \leq I}
} \, \| \, a^j - a^i \, \|_2 \, \right],
\]
which provides an alternative min-max representation of the piecewise affine function $\theta$ with affine pieces $\theta_i$;
this is distinct from the max-min representation of such a piecewise function in \cite{Ovchinnikov02,Scholtes12}.  Third, when
$Q$ is positive semi-definite, $\mbox{dom}(Q,D)$ is known to be convex; in fact, it is equal to the
polyhedron $\left[ \begin{array}{cc}
Q & -D^{\, T} \\
D & 0
\end{array} \right] \left( \mathbb{R}^m \times \mathbb{R}^k_+ \right) - \left( \{ 0 \} \times \mathbb{R}^k_+ \right)$, which is the set
of pairs $(q,b)$ for which the linear constraints (without the complementarity condition) of the KKT system (\ref{eq:KKT QP}) are feasible.
Fourth, there are copositive, nonconvex QPs for which $\mbox{dom}(Q,D)$ is convex.  For instance, for a symmetric, (entry-wise)
positive matrix $Q$ and an identity matrix $D$, $\mbox{dom}(Q,D)$ is equal to the entire space $\mathbb{R}^{m+k}$.  More generally,
if $D_{\infty} = \{ 0 \}$, then for any symmetric matrix $Q$,
$\mbox{dom}(Q,D) = \mathbb{R}^m \times \left[ \, D\mathbb{R}^m - \mathbb{R}^k_+ \, \right]$ is a convex polyhedron.
Fifth, the convexity of ${\cal S}$ in Proposition~\ref{pr:Maher PQ} is a reasonable assumption
because convex, and thus dc, functions are defined only on convex sets.  Hence, the convexity requirement of $\mbox{dom}(Q,D)$
in Theorem~\ref{th:full dc} is needed for one to speak about the dc property of $\mbox{qp}_{\rm opt}(q,b)$.
Admittedly, the dc representation in Theorem~\ref{th:full dc} is fairly complex, due to the possibly exponentially many quadratic pieces
of the value function (cf.\ the proof of Proposition~\ref{pr:copositive results}).  This begs the question of whether a much simpler representation
exists when the matrix $Q$ is positive semi-definite (cf.\ the rather straightforward representation (\ref{eq:pd QP}) in the positive definite case).
There is presently no resolution to this question.

\section{Univariate Folded Concave Functions} \label{sec:folded concave}

The family of univariate folded concave functions was introduced \cite{FanXueZou14} in the literature of sparsity representation
as approximations of the univariate non-zero count function $\ell_0(t) \triangleq \left\{ \begin{array}{ll}
1 & \mbox{if $t \neq 0$} \\
0 & \mbox{otherwise.}
\end{array} \right.$  Formally, such a function is given by $\theta(t) \triangleq f( | t |)$, where $f$ is a continuous, univariate
concave function defined on $\mathbb{R}_+$.
Since we take the domain of $f$ to be the closed interval $[ 0, \infty )$, the composition
property of dc functions \cite[Theorem~II, page~708]{Hartman59} is not applicable to directly deduce that $\theta$ is dc.
We formally state and prove the following
result that is a unification of all the special cases discussed in \cite{AhnPangXin17}.  Proposition~6 in \cite{LeThiPhamVo15}
gives a different dc decomposition of such a folded concave function $\theta(t)$.

\begin{proposition} \label{pr:folded concave} \rm
Let $f$ be a (continuous) univariate concave function defined on $\mathbb{R}_+$.  The composite function $\theta(t) \triangleq f( | t |)$
is dc on $\mathbb{R}$ if and only if
$f^{\, \prime}(0;+)$ exists and is finite, where
\begin{equation} \label{eq:dd folded concave}
f^{\, \prime}(0;+) \, \triangleq \, \displaystyle{
\lim_{\tau \downarrow 0}
} \, \displaystyle{
\frac{f(\tau) - f(0)}{\tau}
} \, = \, -\displaystyle{
\lim_{\tau \uparrow 0}
} \, \displaystyle{
\frac{f(-\tau) - f(0)}{\tau}
} \, .
\end{equation}
\end{proposition}

\begin{figure}[htbp]
\vspace{0cm}\centering{\vspace{0cm}
\includegraphics[scale = .65]{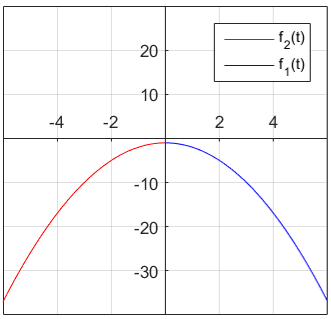}
\includegraphics[scale = .65]{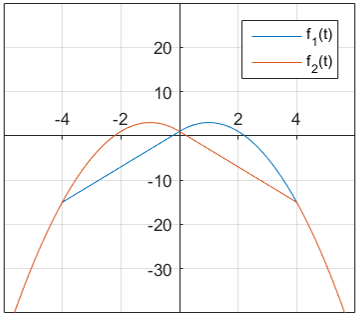}
\includegraphics[scale = .65]{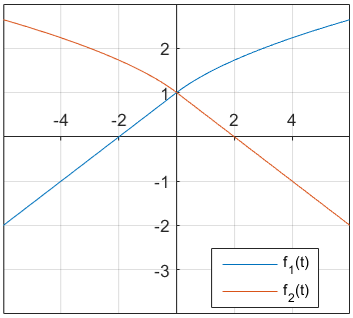}
} \caption{$g(t)= \max\left( f_1(t),f_2(t) \right)$ \\ [5pt]
Case (a): $f(|t|)=-|t|^2-1$;
case (b)1: $f(|t|)=- 2(|t|-1)^2+3$;
case (b)2: $f(|t|)=\sqrt{(|t|+1)}$}\label{fig:folded concavity}
\end{figure}

\noindent {\bf Proof.} Since a dc function must be 
directionally differentiable, a property inherited from a convex
function, it suffices to prove the sufficiency claim of the result.  Since $f$ is concave on $\mathbb{R}_+$, it follows that
\begin{equation} \label{eq:cvity and dd}
f(t) \, \leq \, f(0) + f^{\, \prime}(0;+) \, t, \epc \forall \, t \, \geq \, 0. 
\end{equation}
The proof is divided into 2 cases: (a) $f^{\, \prime}(0;+) \leq 0$, or (b) $f^{\, \prime}(0;+) > 0$.
(See Figure~\ref{fig:folded concavity} for illustration.) In case (a), it follows that
$t = 0$ is a maximum of the function $\theta$ on the interval $( -\infty,\infty )$ by (\ref{eq:cvity and dd}).
Further, we claim that the function $\theta$ is concave on the real line in this case.  To prove the claim, it suffices to
show that if $t_1 > 0 > t_2$, then the secant, denoted $S$, joining the two points $( t_1, \theta(t_1) )$ and $( t_2, \theta(t_2) )$, where
$\theta(t_1) = f(t_1)$ and $\theta(t_2) = f(-t_2)$, on the curve
of $\theta(t)$ is below the curve $\theta(t)$ itself for $t$ in the interval $( t_2, t_1 )$.  This can be argued as follows.  The
secant $S$ can be divided into two sub-secants, one, denoted $S_1$, starting at the end point $( t_1, \theta(t_1) )$, and ending at
$\left( 0, \theta(t_2) - \displaystyle{
\frac{\theta(t_1) - \theta(t_2)}{t_1 - t_2}
} \, t_2 \right)$; and the other, denoted $S_2$, starting at the latter end point and ending at $( t_2, \theta(t_2) )$.
It is not difficult to see that the sub-secant $S_1$ lies below the line segment joining $( t_1, \theta(t_1) )$ and $(0,\theta(0))$,
which in turn lies below the curve of $\theta(t)$ for $t \in ( t_1, 0)$ by concavity of $f$; furthermore, for the same token,
the sub-secant $S_2$ lies below the line segment joining $(0,\theta(0))$ to $( t_2, \theta(t_2) )$, which in turn lies below the curve
of $\theta(t)$ for $t \in ( 0, t_2 )$.  This establishes the concavity of $\theta$ in case (a).

\gap

Consider case (b); i.e., suppose $f^{\, \prime}(0;+) > 0$.   Consider the half-line $t \mapsto f(0) + f^{\, \prime}(0;+) t$ emanating from
the point $(0,f(0))$ for $t < 0$.  Let $t_-^* < 0$ be the right-most $t < 0$ such that this line meets the curve $\theta(t) = f(-t)$
to the left of the origin.  If this does
not happen, we let $t_-^* = -\infty$.  Define the function:
\[ \begin{array}{ll}
f_1(t) \, \triangleq \, \left\{ \begin{array}{ll}
f(t) & \mbox{if $t \, \geq \, 0$} \\ [5pt]
f(0) + f^{\, \prime}(0;+) t & \mbox{if $t \, \in \, [ \, t_-^*,0 \, ]$} \\ [5pt]
f(-t) & \mbox{if $t \, < \, t_-^*$},
\end{array} \right\} & \mbox{if $t_-^* > -\infty$} \\ [0.35in]
f_1(t) \, \triangleq \, \left\{ \begin{array}{ll}
f(t) & \mbox{if $t \, \geq \, 0$} \\ [5pt]
f(0) + f^{\, \prime}(0;+) t & \mbox{if $t \, \in \, ( \, -\infty, 0 \, ]$}
\end{array} \right\} & \mbox{if $t_-^* = -\infty$}
\end{array} \]
We claim that this function is concave on $\mathbb{R}$ and $f_1(t) \leq f(-t)$ for $t \in ( \, t_-^*,0 \, ]$.
Indeed, since $f^{\, \prime}(0;+) > 0$, it follows that
$f(-t)  > f(0) + f^{\, \prime}(0;+) \, t$ for all $t \in ( t_-^*,0 )$ by the definition of
$t_-^*$.  The concavity of $f_1(t)$ can be proved in a way similar to the above proof of case (a), by considering
sub-segments.  Details are omitted.   Similarly, define $t_+^* > 0$ as the left-most $t > 0$ such that half-line
$t \mapsto f(0) - f^{\, \prime}(0;+) \, t$ meets the curve $\theta(t) = f(t)$ to the right of the origin and
let $t_+^* = \infty$ if this does not happen.  Define the function:
\[ \begin{array}{ll}
f_2(t) \, \triangleq \, \left\{ \begin{array}{ll}
f(-t) & \mbox{if $t \, \leq \, 0$} \\ [5pt]
f(0) - f^{\, \prime}(0;+) t & \mbox{if $t \, \in \, [ \, 0, t_+^* \, ]$} \\ [5pt]
f(t) & \mbox{if $t \, > \, t_+^*$},
\end{array} \right\} & \mbox{if $t_+^* < \infty$} \\ [0.35in]
f_2(t) \, \triangleq \, \left\{ \begin{array}{ll}
f(-t) & \mbox{if $t \, \leq \, 0$} \\ [5pt]
f(0) - f^{\, \prime}(0;+) t & \mbox{if $t \, \in \, [ \, 0, \infty \, )$}
\end{array} \right\} & \mbox{if $t_+^* = \infty$}
\end{array} \]
We can similarly show that $f_2(t)$ is concave on $\mathbb{R}$ and $f_2(t) \leq f(t)$ for $t \in [ \, 0, t_+^* \, )$.
Now define $g(t) \triangleq \max\left( f_1(t),f_2(t) \right)$ for all $t \in \mathbb{R}$.  As the pointwise maximum
of two concave
functions, $g$ is dc.  It remains to show that $g(t) = \theta(t)$ for all $t \in \mathbb{R}$.  This can be divided
into 2 cases: $t \geq 0$ and $t \leq 0$.  In each case, the above established
properties of the two functions $f_1$ and $f_2$ can be applied to complete the proof.  \hfill $\Box$

\gap

%

\gap

{\bf Remarks.} Being not dc, the univariate function $\theta(t) \triangleq \sqrt{| \, t \, |}$ provides a counter-example
to illustrate the important role of the existence of the limit (\ref{eq:dd folded concave}).
Another relevant remark is that the following fact is known \cite[page~707]{Hartman59}: ``if
${\cal D}$  is a (bounded or unbounded)
interval, then the univariate function $f$ is dc on ${\cal D}$ if and only if $f$ has left and right
derivatives (where these are meaningful) and these derivatives are of
bounded variation on every closed bounded interval interior to ${\cal D}$''.  We have not applied this fact to
prove Proposition~\ref{pr:folded concave} because our proof provides a simple construction of the dc representation
of the function $\theta$ in terms of the function $f$.  \hfill $\Box$

\gap

{\bf Acknowledgement.}  The second author gratefully acknowledges the discussion with Professors Le Thi Hoi An and Pham Dinh Tao
in the early stage of this work during his visit to the Universit\'e de Lorraine, Metz in June 2016.
The three authors acknowledge their fruitful discussion with Professor Defeng Sun at the National
University of Singapore during his visit to the University of Southern California.  They are also grateful to Professor Marc Teboulle
for drawing their attention to the references \cite{BenTalTeboulle07,BenTalTeboulle87,BenTalTeboulle86} that introduce and revisit
the OCE. The authors are also grateful to Dr. Ying Cui for bringing to our attention the reference \cite{Klatte85} and Example 4.1 therein. The constructive comments of two referees are also gratefully acknowledged.
In particular, the authors are particularly grateful to a referee
who has been very patient with their repeated misunderstanding of the work \cite{Wozabal12} that is now correctly summarized at the end of Subsection~\ref{subsec:CVaR and VaR}.

\end{document}